\def\a{\alpha}

\def\cd{\cdot}
\def\d{\delta}

\def\e{\epsilon}

\def\G{\Gamma}

\def\l{\ell}
\def\lb\{{\left\{}
\def\la{\lambda}
\def\La{\Lambda}
\def\lla{\longleftarrow}
\def\lm{\limits}
\def\lra{\longrightarrow}
\def\dllra{\Longleftrightarrow}
\def\llra{\longleftrightarrow}
\def\n{\nabla}
\def\ngth{\negthickspace}
\def\ola{\overleftarrow}
\def\Om{\Omega}
\def\om{\omega}
\def\op{\oplus}
\def\oper{\operatorname}
\def\oplm{\operatornamewithlimits}
\def\ora{\overrightarrow}
\def\ov{\overline}
\def\ova{\overarrow}
\def\ox{\otimes}
\def\p{\partial}
\def\rb\}{\right\}}
\def\s{\sigma}
\def\sbq{\subseteq}
\def\spq{\supseteq}
\def\sqp{\sqsupset}
\def\supth{{\text{th}}}
\def\T{\Theta}
\def\th{\theta}
\def\tl{\tilde}
\def\thra{\twoheadrightarrow}
\def\un{\underline}
\def\ups{\upsilon}
\def\vp{\varphi}
\def\wh{\widehat}
\def\wt{\widetilde}
\def\x{\times}
\def\z{\zeta}
\def\({\left(}
\def\){\right)}
\def\<{\left<}
\def\>{\right>}

\def\SA{\mathcal A}
\def\SB{\mathcal B}
\def\SC{\mathcal C}
\def\SD{\mathcal D}
\def\SE{\mathcal E}
\def\SF{\mathcal F}
\def\SG{\mathcal G}
\def\SH{\mathcal H}
\def\SI{\mathcal I}
\def\SJ{\mathcal J}
\def\SK{\mathcal K}
\def\SL{\mathcal L}
\def\SM{\mathcal M}
\def\SN{\mathcal N}
\def\SO{\mathcal O}
\def\SP{\mathcal P}
\def\SQ{\mathcal Q}
\def\SR{\mathcal R}
\def\SS{\mathcal S}
\def\ST{\mathcal T}
\def\SU{\mathcal U}
\def\SV{\mathcal V}
\def\SW{\mathcal W}
\def\SX{\mathcal X}
\def\SY{\mathcal Y}
\def\SZ{\mathcal Z}

\newcommand{\BA}{\ensuremath{\mathbf A}}
\newcommand{\BB}{\ensuremath{\mathbf B}}
\newcommand{\BC}{\ensuremath{\mathbf C}}
\newcommand{\BD}{\ensuremath{\mathbf D}}
\newcommand{\BE}{\ensuremath{\mathbf E}}
\newcommand{\BF}{\ensuremath{\mathbf F}}
\newcommand{\BG}{\ensuremath{\mathbf G}}
\newcommand{\BH}{\ensuremath{\mathbf H}}
\newcommand{\BI}{\ensuremath{\mathbf I}}
\newcommand{\BJ}{\ensuremath{\mathbf J}}
\newcommand{\BK}{\ensuremath{\mathbb K}}
\newcommand{\BL}{\ensuremath{\mathbf L}}
\newcommand{\BM}{\ensuremath{\mathbf M}}
\newcommand{\BN}{\ensuremath{\mathbf N}}
\newcommand{\BO}{\ensuremath{\mathbf O}}
\newcommand{\BP}{\ensuremath{\mathbf P}}
\newcommand{\BQ}{\ensuremath{\mathbb Q}}
\newcommand{\BR}{\ensuremath{\mathbf R}}
\newcommand{\BS}{\ensuremath{\mathbf S}}
\newcommand{\BT}{\ensuremath{\mathbf T}}
\newcommand{\BU}{\ensuremath{\mathbf U}}
\newcommand{\BV}{\ensuremath{\mathbf V}}
\newcommand{\BW}{\ensuremath{\mathbf W}}
\newcommand{\BX}{\ensuremath{\mathbf X}}
\newcommand{\BY}{\ensuremath{\mathbf Y}}
\newcommand{\BZ}{\ensuremath{\mathbb Z}}

\def\bba{{\mathbb A}}
\def\bbb{{\mathbb B}}
\def\bbc{{\mathbb C}}
\def\bbd{{\mathbb D}}
\def\bbe{{\mathbb E}}
\def\bbf{{\mathbb F}}
\def\bbg{{\mathbb G}}
\def\bbh{{\mathbb H}}
\def\bbi{{\mathbb I}}
\def\bbj{{\mathbb J}}
\def\bbk{{\mathbb K}}
\def\bbl{{\mathbb L}}
\def\bbm{{\mathbb M}}
\def\bbn{{\mathbb N}}
\def\bbo{{\mathbb O}}
\def\bbp{{\mathbb P}}
\def\bbq{{\mathbb Q}}
\def\bbr{{\mathbb R}}
\def\bbs{{\mathbb S}}
\def\bbt{{\mathbb T}}
\def\bbu{{\mathbb U}}
\def\bbv{{\mathbb V}}
\def\bbw{{\mathbb W}}
\def\bbx{{\mathbb X}}
\def\bby{{\mathbb Y}}
\def\bbz{{\mathbb Z}}

\documentclass[12pt]{amsart}
\usepackage{amsmath,amsthm}
\usepackage{amssymb,pb-diagram,enumerate}
\usepackage{float,epsfig}
\usepackage{psfrag}
\usepackage{graphicx}

\newtheorem{thm}{Theorem}[section]
\newtheorem{lem}[thm]{Lemma}
\newtheorem{prop}[thm]{Proposition}
\newtheorem{cor}[thm]{Corollary}
\newtheorem{conj}{Conjecture}\renewcommand{\theconj}{}

\theoremstyle{definition}
\newtheorem{defn}[thm]{Definition}
\newtheorem{rem}[thm]{Remark}

{\catcode`@=11\global\let\c@equation=\c@thm}
\renewcommand{\theequation}{\thethm}

{\catcode`@=11\global\let\c@figure=\c@thm}
\renewcommand{\thefigure}{\thethm}

\newcommand{\1}{\ensuremath{\pi_1}}
\newcommand{\2}[1]{\ensuremath{^{(#1)}}}
\newcommand{\nn}{\ensuremath{^{(n+1)}}}
\newcommand{\ff}{\ensuremath{\SF_{(n)}/\SF_{(n.5)}}}
\newcommand{\zz}{\ensuremath{\SG_{(n+2)}/\SG_{(n+2.5)}}}
\newcommand{\sr}{\ensuremath{S^3\backslash R}}
\newcommand{\sk}{\ensuremath{S^3\backslash K}}
\newcommand{\ns}{$n$-solvable}
\newcommand{\nss}{$(n.5)$-solvable}

\newcommand{\arf}{\operatorname{Arf}}
\newcommand{\Hom}{\operatorname{Hom}}
\newcommand{\id}{\operatorname{id}}
\newcommand{\image}{\operatorname{image}}
\newcommand{\kernel}{\operatorname{kernel}}
\newcommand{\rank}{\operatorname{rank}}
\newcommand{\sign}{\operatorname{sign}}
\newcommand{\ind}{\operatorname{index}}
\newcommand{\tr}{\operatorname{trace}}
\newcommand{\vol}{\operatorname{volume}}

\begin{document}

\title{Knot concordance and von Neumann {\boldmath $\rho$}-invariants}
\author{Tim D. Cochran}
\address{Rice University, Houston, Texas, 77005-1892}
\email{cochran@math.rice.edu}
\author{Peter Teichner}
\address{University of California, Berkeley, CA, 94720-3840}
\email{teichner@math.berkeley.edu}
\thanks{Both authors are partially supported by the National Science Foundation}
\begin{abstract}
We prove the nontriviality at all levels of the filtration of the classical topological knot concordance group $\SC$
$$
\cdots \subseteq \SF_{n} \subseteq \cdots \subseteq
\SF_1 \subseteq\SF_{0} \subseteq \SC.
$$
defined in \cite{COT}. This filtration is significant because not only is it strongly connected to Whitney tower constructions of Casson and Freedman, but all previously-known concordance invariants are related to the first few terms in the filtration. In \cite{COT} we proved nontriviality at the first new level $n=3$ by using von Neumann $\rho$-invariants of the $3$-manifolds obtained by surgery on the knots. For larger $n$ we use the Cheeger-Gromov estimate for such $\rho$-invariants, as well as some rather involved algebraic arguments using our noncommutative Blanchfield forms.
In addition, we consider a closely related filtration, $\{\SG_n\}$, of $\SC$ defined in terms of  Gropes in the $4$-ball. We show that this filtration is also non-trivial for all $n>2$.
\end{abstract}
\maketitle

\section{Introduction}
\begin{figure}[htbp]
\setlength{\unitlength}{1pt}
\begin{picture}(141,111)
\put(53,98){$K$}
\put(10,10){\includegraphics{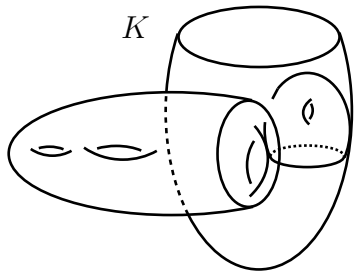}}
\end{picture}
\caption{A  Grope of height $2$.} \label{Grope}
\end{figure}

In \cite{COT}  Kent Orr and the authors exhibit a
new geometrically defined filtration of the classical knot concordance group 
$\SC$
$$
\cdots \subseteq \SF_{n} \subseteq \cdots \subseteq
\SF_1\subseteq \SF_{0.5} \subseteq \SF_{0} \subseteq \SC.
$$ 
For each positive half-integer $n$, the subgroup $\SF_{n}$ consists of all {\em
$(n)$-solvable} knots. $\SF_0$ consists of the knots with Arf-invariant zero, and $\SF_{0.5}$ consists of the algebraically slice knots. In general, $(n)$-solvability is defined using intersection forms in $(n)$-solvable covering spaces of certain $4$-manifolds (the definition is reviewed in Section~\ref{sec:not solvable}). Here, we also consider a closely related filtration that is geometrically more intuitive, namely a filtration defined using \emph{Gropes}. Denote by $\SG_n$ the subgroup of $\SC$ consisting of all knots that bound a Grope of height~$n$ in $D^4$. The precise definitions are reviewed in Section~\ref{sec:gropes} but a grope of height $2$ is shown in Figure~\ref{Grope}. Note that we use exclusively {\em symmetric} gropes, corresponding to the derived series of a group (see \cite{T} for other notions of gropes and their applications).

The filtrations of $\SC$ by solvability and by Gropes are not the same but the following result from \cite{COT} shows they are closely correlated and so the reader may safely think of the Grope filtration as a more intuitive approximation to the filtration by solvability.

\begin{prop}\label{oldtheorem} \cite[Theorem 8.11]{COT} If a knot $K$ in $S^3$ bounds a Grope of height~$(n+2)$ in $D^4$ then $K$ is $(n)$-solvable.  Thus, $\SG_{n+2} \subseteq \SF_n$. 
\end{prop}

The main result of this paper, and in some sense the ultimate justification for the entire theory of \cite{COT} is  

\begin{thm} \label{thm:main}
For any $n \in \bbn_0$, the quotient groups $\SF_{n}/\SF_{n.5}$ contain elements of infinite order. Similarly, the groups $\SG_{n+2}/\SG_{n+2.5}$ contain elements of infinite order. In fact, the groups $\SG_{n+2}$ contain knots that have infinite order modulo $\SF_{n.5}$.
\end{thm}
For $n=0$, this is detected by the Seifert form obstructions \cite[Remark 1.13]{COT}, for $n = 1$ one can use Casson-Gordon invariants (that vanish on $\SF_{1.5}$ by \cite[Theorem 9.11]{COT}) and for $n=2$ this is the main result of \cite[Theorem 6.4]{COT}. The results for $n>2$ are new.

The proof of this theorem turns out to be technically quite difficult. The main problem is the obstruction theoretic nature of all known invariants, a problem that already arises in using Casson-Gordon invariants. In a nutshell, there \emph{are} signature invariants that obstruct the existence of a Grope of height $n.5$ {\em extending} a given Grope of height~$n$. But, in order to show that a knot does not bound a Grope of height $n.5$ one then has to show that this signature is nontrivial {\em for all possible} Gropes of height~$n$ bounding the knot. This turns out to be a formidable task. It is resolved by Casson-Gordon and in \cite[Section 6]{COT} by finding knots for which, roughly speaking, there is a {\em unique} Grope of height~$n$, at least as far as the relevant algebra can see. This works for $n=1,2$ but we have been unable to find such knots for $n>2$. In this paper we use analytic tools to resolve the issue. Of course we also use the noncommutative localization techniques of \cite{COT} and ideas from \cite{COT2} to construct the relevant knots, with a little more care necessary to make sure that they lie in $\SG_{n+2}$, rather than just in $\SF_n$. The new analytic methods show that these knots are not \nss.

To see why analysis can play a role, recall that the signature that turns out to be relevant for $n>1$ is in fact a (real valued) von Neumann signature, associated to a certain intersection form on a 4-manifold constructed from a Grope of height~$n$. This 4-manifold has boundary $M_K$, the $0$-surgery on the given knot $K$. The main information about the Grope is encoded by the homomorphism 
\[
\phi:\1M_K \lra \pi
\]
induced by the inclusion of the boundary into the 4-manifold. Since the Grope is not variable, we do not have any information on the group $\pi$, in particular, we don't know any interesting representations, except for the canonical one on $\l^2(\pi)$. That's why von Neumann algebras enter the story.

By the von Neumann index theorem, the difference between the von Neumann signature and the usual (untwisted) signature of this $4$-manifold is equal to the invariant $\rho(M_K,\phi)$ of the boundary. This von Neumann $\rho$-invariant is the difference between the von Neumann $\eta$-invariant and the untwisted $\eta$-invariant. It is a real-valued topological invariant, depending only on the covering determined by $\phi$ and will be explained in detail in Section~\ref{sec:CG}.

Thus the main technical question of how to control the von Neumann signatures of the $4$-manifolds associated to all possible Gropes of height~$n$ with boundary $K$, translates into the question of how to control the invariants $\rho(M_K,\phi)$ for all possible $\phi$. The parameter $n$ enters in this reformulation since $\pi$ is really the quotient of the 4-manifold group by $(n+1)$-fold commutators. Analysis enters prominently because of the following estimate of Cheeger and Gromov \cite{CG2} for the von Neumann $\rho$-invariants
\begin{equation}\label{CG}
\exists \quad C_M>0 \text{ such that } |\rho(M,\phi)| < C_M \quad \forall \phi
\end{equation}
where $M$ is a closed oriented manifold of odd dimension. Thus, for any fixed $3$-manifold, the set of all possible $\rho$-invariants is bounded in absolute value. In addition to the original paper, Ramachandran's paper \cite{R} is a good source to learn about the estimate in more generality (see our Section~\ref{sec:CG}).

\subsection*{Outline of the proof of Theorem~\ref{thm:main}}  We start with the $0$-surgery $M$ on a fibred genus~2 ribbon knot $R$ (which bounds an embedded disk in $B^4$ and is thus $(n)$-solvable for all $n$). We prove the existence of a certain collection of circles $\eta_1,\dots, \eta_m$ in the $n$-th derived subgroup of $\1M$ that forms a trivial link in $S^3$. We then modify $R$ (we call this a \emph{genetic infection}) using a certain auxiliary knot $J$ (called the \emph{infection knot}) along the circles  $\eta_1,\dots, \eta_m$ (called \emph{axes}).

\begin{figure}[htbp]
\setlength{\unitlength}{1pt}
\begin{picture}(265,75)
\put(5,41){$\eta_1$}
\put(53,41){$\dots$}
\put(123,41){$\eta_m$}
\put(182,36){$J$}
\put(238,37){$J$}
\put(30,6){$R$}
\put(85,6){$R$}
\put(180,6){$R(\eta_1,\dots,\eta_m,J)$}
\put(207,41){$\dots$}
\put(20,20){\includegraphics{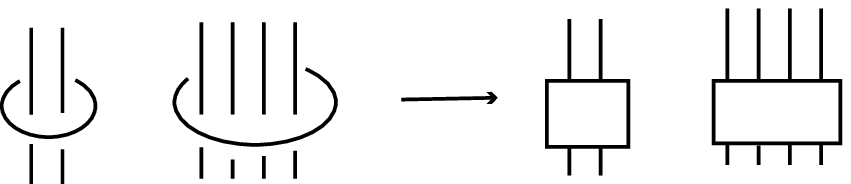}}
\end{picture}
\caption{Genetic infection of $R$ by
$J$ along $\eta_i$}\label{infection}
\end{figure}

\begin{figure}[htbp]
\setlength{\unitlength}{1pt}
\begin{picture}(268,103)
\put(0,70){$\eta_1$}
\put(107,70){$\eta_2$}
\put(178,55){$J$}
\put(245,56){$J$}
\put(60,5){$R$}
\put(178,5){$K=R(\eta_1,\eta_2,J)$}
\put(0,20){\includegraphics{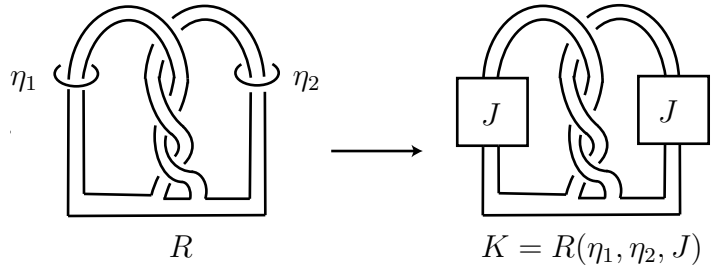}}
\end{picture}
\caption{An example for $n=1$}\label{example}
\end{figure}

This language, introduced in \cite{COT2}, just expresses a standard construction, where the $\eta_i$ bound disjointly embedded disks in $S^3$ and the knot $J$ is tied into all strands of $R$ passing through one of these disks as illustrated in general in Figure~\ref{infection} and in a particular case in Figure~\ref{example}, where $n=1$, $m=2$.

With a good choice of an infection knot $J$, we can make sure that all resulting knots $K=R(\eta_1,\dots, \eta_m,J)$ bound Gropes of height $(n+2)$, see Theorem~\ref{thm:examples}. 
\begin{figure}[htbp]
\setlength{\unitlength}{1pt}
\begin{picture}(260,150)
\put(10,121){$U$}
\put(247,74){$J$}
\put(203,78){$\cong$}
\put(20,20){\includegraphics{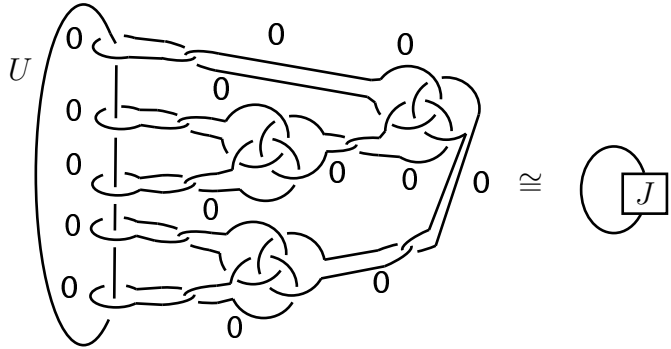}}
\end{picture}
\caption{The knot $J$ as surgery on an unknot $U$}\label{knotJ}
\end{figure}

Figure~\ref{knotJ} shows, implicitly, our choice of the infection knot $J$ (in fact, we use the connected sum of many copies of $J$, in order to get a large signature, see below). The surgery point of view that defines $J$ turns out to have enormous flexibility for our purposes. In fact, the picture for $J$ will be redrawn 6 times in this paper, depending exactly on the aspect of knot theory that we study. The first instance is Figure~\ref{JClasper} which explains the clasper origin of our knot.

If the infinite cyclic von Neumann signature of the infection knot is bigger than $C_M$ (the ``infection is strong enough'') then the Cheeger-Gromov estimate (\ref{CG}) shows that $K$ is not \nss, provided that for all (relevant) homomorphisms $\phi$, one of the axes $\eta_i$ is mapped nontrivially. The precise statement is given in Theorem~\ref{thm:infection}.

The last condition on the axes $\eta_i$ turns out to be more subtle than expected, mainly because $\1M^{(n)}/\1M^{(n+1)}$ is not finitely generated for $n>1$. If it were then any generating set would have the desired property. Fortunately, it turns out that we can use the higher order Blanchfield forms from \cite{COT} to carefully construct axes $\eta_i$ with the desired property, see Theorem~\ref{thm:injectivity}.

The main remaining open questions about our filtrations is
\begin{conj} For any $n \in \bbn_0$, the quotient groups $\SF_{n}/\SF_{n.5}$ have infinite rank. Moreover, for $n>0$, the groups $\SF_{n.5}/\SF_{n+1}$ are nontrivial.
\end{conj}

\noindent The first part of the conjecture is true for $n\leq 2$:
$n = 0$ is verified using twist knots and the Seifert form obstructions, for $n = 1$ this
can be established by using examples due to Casson-Gordon,
and $n=2$ is the main result of \cite{COT2}. 

Our paper is organized as follows: In Section~\ref{sec:CG}, we give a survey of the analytical results surrounding the von Neumann $\rho$-invariant; no originality is claimed. In Theorem~\ref{thm:examples} we construct a large class of knots that bound Gropes of height~$(n+2)$ in the 4-ball. Theorems~\ref{thm:infection} and \ref{thm:injectivity} together imply that many of these knots are not $(n.5)$-solvable. Therefore, our main result Theorem~\ref{thm:main} follows. Section~\ref{sec:noncommutative} reviews from \cite{COT}  some algebraic topological results, higher-order Alexander modules and their Blanchfield forms, used for the main technical Theorem~\ref{thm:paininbutt} that is at the heart of the proof of Theorem~\ref{thm:injectivity}. 

This paper uses many notions and ideas from \cite{COT}. However, we have made a conscious effort to include here all relevant definitions, and we stated all results needed explicitly, with precise references. We hope that, as a consequence, most readers will be able to enjoy this paper independently of \cite{COT}. Those who want to see the proofs of the stated results will find them easily in \cite{COT}.

\vspace{0.5 cm}
\noindent {\em Acknowledgements}: 
It is a pleasure to thank Bruce Driver and Lance Small for discussions on some of the analytic, respectively algebraic, aspects of the paper, and Kent Orr for several useful conversations.

\newpage

\section{Von Neumann $\rho$-invariants} \label{sec:CG}

For the convenience of the reader, this section gives a short survey of the more analytic aspects of the invariant that is used in the rest of the paper. Throughout the section we work with the signature operator on 4-dimensional manifolds (and their boundary), even though everything works as well on $4k$-manifolds. One needs to replace $p_1/3$ by Hirzebruch's L-polynomial (and the 1-forms below by $(2k-1)$ forms on the boundary). In fact, most of our discussion applies to all Dirac type operators in any dimension, instead of just the signature operator. 

\subsection*{The signature theorem for manifolds with boundary}
Let $W$ be a compact oriented Riemannian $4$-manifold with boundary $M$ and assume that the metric is a product near the boundary (or at least that the first two normal derivatives of the metric vanish on $M$). The Atiyah-Patodi-Singer index theorem for the case of the signature operator  implies the following signature theorem \cite[Thm.4.14]{APS}
\begin{equation} \label{APS}
\s(W) =   \int_W p_1(W)/3 - \eta(M) 
\end{equation}
where $\s(W)$ is the signature of the intersection form on $H_2(W)$, $p_1(W)$ is the first Pontrjagin form of the tangent bundle (which depends on the metric), and $\eta(M)$ is a spectral invariant of the boundary. It is the value at $s=0$ of the $\eta$-function
\begin{equation}\label{eta}
\eta(s)= \sum_{\la\neq 0} (\sign \la) |\la|^{-s} .
\end{equation}
Here $\la$ runs through the nonzero Eigenvalues of the signature operator $D$ on $M$. More precisely, $D$ is the self-adjoint operator on \emph{even} differential forms on $M$ defined by $\pm(*d-d*)$ where the Hodge $*$-operator depends on the metric. The $\eta$-function is defined by analytic continuation and it turns out that it is holomorphic for $Re(s)>-\frac 12$. 
In fact, one can get the explicit formula, used below,
\begin{equation} \label{gamma}
\eta(2s) = \frac{1}{\G(s+1/2)} \int_0^\infty t^{s-1/2} \tr(De^{-tD^2}) dt
\end{equation}
where the trace class operator $De^{-tD^2}$ is defined by functional calculus and the Gamma function is given by
\[
\G(s)=  \int_0^\infty t^{s-1} e^{-t} dt
\]
Using Hodge theory, \cite[Prop.4.20]{APS} show that for the purpose of calculating the $\eta$-function, one may restrict to the operator $d*$ acting on the space $d\Om^1(M)$. As pointed out by \cite{APS} this translation leads to the following suggestive interpretation of $\eta(M)$. Define a quadratic form $Q$ by
\begin{equation}\label{Q}
Q(\alpha) := \int_M \a \wedge d\a, \quad \a\in\Om^1(M)
\end{equation}
and observe that $Q$ has radical $\ker(d)$ and hence gives a form on $d\Om^1(M)$. Moreover, for Eigenvectors $d\a$ of $d*$, the corresponding Eigenvalue has opposite sign as $Q(\a)$. Hence one can formally interpret the correction term $-\eta(M)$ in (\ref{APS}) as the ``signature'' of the quadratic form $-Q$. Note that $Q$ does \emph{not} depend on the metric but since $Q$ is defined on the infinite dimensional space, the metric is used to give the proper meaning to its signature.

Even though the above signature theorem~(\ref{APS}) is what we need in this paper, it might be good to remind the reader of the relation to the index of the signature operator on $W$. In fact, the index theorem for the signature operator $D_W$ on $W$ reads as follows \cite[4.3]{APS}
\begin{equation} \label{index}
\ind(D_W) =  \int_W  p_1(W)/3 -(h+ \eta(0))/2
\end{equation}
where $h$ is the dimension of the space of harmonic forms on $M$ and $\eta(0)$ is the value at $s=0$ of the $\eta$-function for the signature operator on \emph{all} forms on $M$. Since this operator preserves the parity of forms and commutes with the Hodge $*$-operator, the above $\eta$-function is just \emph{twice} the $\eta$-function (on even forms) considered in the signature theorem. This explains the disappearance of the factor $1/2$ from (\ref{index}). Moreover, by considering $L^2$-solutions on $W$ with an infinite cylinder attached, \cite[4.8]{APS} show that
\[
\ind(D_W) = \s(W) - h
\]
which explains how (\ref{APS}) is derived from (\ref{index}).

In \cite{CG1} Cheeger and Gromov were looking for geometric conditions under which the integral term in the signature theorem
is a topological invariant. They proved the following beautiful result: Assume the open $4$-manifold $X$ admits a \emph{complete} metric with finite volume whose sectional curvature satisfy
\[
|K(X)| \leq 1 
\]
Assume furthermore that $X$ has \emph{bounded covering geometry} in the sense that there is a normal covering $\tilde X$ (with covering group $\G$) such that the injectivity radius of the pull back metric on $\tilde X$ is $\geq 1$ (this condition gets weaker with the covering getting larger). 
Then \cite[Thm.6.1]{CG1} says that one has a topological invariant (in fact, a proper homotopy invariant)
\begin{equation}
\int_X p_1(X)/3 = \s_\G(X)
\end{equation}
where $\s_\G(X)$ is the \emph{von Neumann-}, or $L^2$-signature of the $\G$-cover $\tl X$ to which we turn in the next section. We point out that the main step in the proof of this theorem is an estimate for the $\eta$-invariant of $3$-manifolds (with similarly bounded geometry), compare Theorem~\ref{thm:CG} below. It will imply what we called the \emph{Cheeger-Gromov estimate} in our introduction.

\subsection*{The {\boldmath $L^2$}-signature theorem}
Returning to a compact $4$-manifold $W$ with product metric near the boundary, we can study twisted signatures, given by bundles with connection over $W$. If the bundle is flat then these signatures have again a homological interpretation (and the integral term in the signature theorem is unchanged). A flat bundle is given by a representation of the fundamental group $\1W$. However, in the application we have in mind, there are no preferred such representations mainly because $\1W$ is an unknown group. All we will be given is a homomorphism $\1W \to \G$ where $\G$ is usually a solvable group. Fortunately, there is a highly nontrivial \emph{canonical} representation of any group, namely on $\ell^2(\G)$. It turns out that one can twist the signature operator with this representation and then calculate its index using von Neumann's $\G$-dimension. We used these real numbers in \cite{COT} to prove our main results and we gave a survey in section~5, similar to the current one, including all relevant references. Since then, L\"uck and Schick wrote \cite{LS} where they prove that all known definitions of von Neumann signatures agree.

There is an $L^2$-signature theorem, analogous to (\ref{APS}) which can again be derived from an $L^2$-index theorem, proven for general Dirac type operators in \cite{R}. A more direct argument for signatures is given in \cite[Thm.3.10]{LS} where the authors use Vaillant's thesis to translate the right hand side of (\ref{L2 APS}) below to the $L^2$-signature of the intersection pairing on $L^2$-harmonic 2-forms on $W$ with an infinite cylinder attached. L\"uck and Schick then translate this signature into a purely homological setting obtaining the following result: Given a compact oriented $4$-manifold $W$ together with a homomorphism $\1W\to\G$ once has
\begin{equation} \label{L2 APS}
\s_\G(W) =   \int_W p_1(W)/3 - \eta_\G(M) 
\end{equation}
Here the left hand side is the von Neumann signature of the intersection form on $H_2(W;\SN\G)$, where $\SN\G$ is the von Neumann algebra of the group $\G$. The von Neumann $\eta$-invariant can also be defined in a straight-forward way, using (\ref{gamma}) rather then (\ref{eta}):
\begin{equation} \label{L2 gamma}
\eta_\G(M) := \frac{1}{\sqrt\pi} \int_0^\infty t^{-1/2} \tr_\G(\tilde De^{-t\tilde D^2}) dt
\end{equation}
Here $\tilde D$ is the signature operator on the even forms of the induced $\G$-cover $\tl M$ of $M$. If $k_t(x,y)$ denotes the smooth kernel of the operator $\tilde De^{-t\tilde D^2}$ then the above $\G$-trace is given by the integral
\[
\int_\SF \tr_xk_t(x,x)dx
\]
where $\SF$ is a fundamental domain for the $\G$-action on $\tl M$.

As before, one can alternatively use the lift of $d*$ to calculate this $\eta_\G$-invariant. 
A key step in the proof of (\ref{L2 APS}) is the following result of \cite[4.10]{CG2}. 
\begin{thm}\label{thm:CG}
There is a constant $C$, depending only on the local geometry of $M$, such that for all $\G$-covers of $M$
\[
|\eta_\G(M)| \leq C \cdot \vol(M)
\]
\end{thm}
This is also proven in \cite[Thm.3.1.1]{R} for general Dirac type operators but there seems to be a problem for $s$ near $t$ in the 4-th line of equation (3.1.10). As Bruce Driver pointed out to us, this problem can be fixed by using an $L^1$-estimate instead of the $L^\infty$ estimate in (3.1.10). For the signature operator $d*$, one can also carefully read pages 23 and 24 of \cite{CG2}, inserting the symbols $*d$ into the decisive definition (4.15) of the von Neumann $\eta$-invariant. Note that Cheeger and Gromov use the operator $*d$ on coexact 1-forms which is conjugate, under Hodge-$*$, to $d*$ on $d\Om^1$ used above.

\subsection*{The von Neumann  {\boldmath $\rho$}-invariant}
Subtracting the expressions (\ref{APS}) and (\ref{L2 APS}) one gets the following equation
\begin{equation}\label{rho}
\s_\G(W) - \s(W) = \eta(M) - \eta_\G(M)
\end{equation}
which allows the following beautiful interpretation: the left hand side is independent of the metric, whereas the right hand side does not depend on the zero bordism $W$ for $M$. As a consequence, the above expression must be a topological invariant of $M$! This argument works as long as, for a given $3$-manifold $M$, one can find a metric (easy) and a zero bordism (also easy, except if it has to be over the group $\G$). In fact, it suffices to find the zero bordism over a group into which $\G$ embeds and this can always be done. Alternatively, one applies (\ref{rho}) to the product $M \times I$, equipped with a metric inducing a path of metrics on $M$. Since the signature terms vanish on the product, one concludes that the right hand side is indeed independent of the metric. 
\begin{defn}\label{def:rho}
Let $M$ be a closed oriented $3$-manifold and fix a homomorphism $\phi:\1M\to \G$. Define the von Neumann $\rho$-invariant
\[
\rho_\G(M,\phi) :=  \eta(M) - \eta_\G(M)
\]
with respect to any metric on $M$. If the group (or the homomorphism) is clear from the context, we suppress it from the notation, like we already did for $\eta_\G$.
\end{defn}
It would be extremely interesting to find a combinatorial interpretation of this von Neumann $\rho$-invariant. Since it does not depend on the choice of a metric, a definition along the lines of the quadratic form $Q$ in (\ref{Q}) might not be out of reach. 
All our calculations of $\rho_\G$ are based on (\ref{rho}) and in fact, by choosing our examples of knot carefully, we manage to reduce the $\G$-signature calculations to the case $\G=\bbz$. There it boils down to an integral, over the circle, of all twisted signatures, one for each $U(1)$ representation of $\bbz$, compare Lemma~\ref{signatureJ}.

The Cheeger-Gromov estimate from Theorem~\ref{thm:CG} now clearly implies the following innocent looking estimate. It will be crucial for our purposes.
\begin{thm}\label{thm:estimate}
For any closed oriented $3$-manifold $M$, there is a constant $C_M$ such that for all groups $\G$ and all homomorphisms $\phi:\1M\to\G$
\[
|\rho_\G(M,\phi)|<C_M.
\]
\end{thm}
We noticed long ago that such an estimate would be extremely helpful for understanding our filtration of the knot concordance group. We made several unsuccessful attempts at proving this estimate by using (\ref{rho}), i.e. the interpretation in terms of signature defects. Our lesson is that estimates are best approached with analytic tools.

\section{Gropes of height $(n+2)$ in $D^4$} \label{sec:gropes}

In this section we review the definition of a \emph{Grope} and then describe, for each positive integer $n$, large families of knots in $S^3$ that bound embedded Gropes of
height $(n+2)$ in $D^4$. In the next section we show that among
these are knots that {\it do not} bound any Grope of
height $(n+2.5)$ in $D^4$.

We now review the definition of a \emph{Grope} (see \cite{FQ} \cite{FT}). More precisely, we shall only define {\em symmetric} gropes (and Gropes) since these are used exclusively. We shall therefore suppress the adjective ``symmetric'', for a survey of other notions and applications, see \cite{T}.

In the following, when we refer to a
``surface'', we mean a compact, connected, oriented surface with one
boundary component. Recall that each connected component of an abstract grope $G$
is built up from a connected \emph{first-stage surface} $G^1$ by
gluing \emph{second-stage} surfaces to each circle in a symplectic
basis $\{a_j,b_j\}$, $1\leqq j\leqq 2g$, for $G^1$. This growth
process continues, so that in general the $k$-stage surface
$(G^k)_{a_j}$ is glued to $a_j$ and $(G^k)_{b_j}$ is glued to
$b_j$ where $\{a_j,b_j\}$ is a symplectic basis of circles for one
of the $(k-1)$-stage surfaces. If a grope has $n$ stages in all then we say that it has \emph{height} $n$. A grope of height $0$ is understood to be merely a circle (no surfaces). The union of all of the circles in symplectic bases for the $n$-stage surfaces is called the set of \emph{tips} of the grope. No surfaces are attached to the tips and they freely generate the fundamental group of the grope.

The reader may have observed that gropes are related to the derived series of a group as follows. Let $A^{(i)}$ denote the {\it $i$-th derived group} of a group $A$,
inductively defined by $A^{(0)}:=A$ and $A^{(i+1)}:=[A^{(i)}, A^{(i)}]$.
A group $A$ is {\it $n$-solvable} if $A^{(n+1)}=1$ ($0$-solvable
corresponds to abelian) and $A$ is {\it solvable} if such a finite $n$
exists. The connection between  gropes and the derived series is then given by the following statement, see \cite[Part II Lemma 2.1]{FT}. A loop $\gamma :S^1\to X$ in a space $X$ extends to a continuous map of a height $n$  grope $G\to X$ if and only if $\gamma$ represents an element of the $n$-th term of the derived series of the fundamental group of $X$. In particular, for $\gamma$ to bound an \emph{embedded} such grope is a strictly stronger statement. In this paper when we say that a circle \emph{bounds a grope} in a space X, we will mean that it bounds an \emph{embedded} copy of the abstract grope. The difference between this geometric condition and its algebraic counterpart is a central underlying theme of this section.

If a grope $G$ is embedded in a 4-dimensional manifold, we usually would like to be able to take arbitrary many disjoint parallel copies, hence we require the following framing condition: A neighborhood of $G$ is diffeomorphic to the product of $\bbr$ with a neighborhood of the easiest embedding of $G$ into $\bbr^3$. Another way of expressing this is to say that the relative Euler number of each surface stage vanishes. This relative invariant is defined because the boundary circle of each surface stage, except the bottom, inherits a framing from its embedding into the previous stage. In \cite{FT},  neighborhoods of such framed gropes in $4$-manifolds were called {\em Gropes} and we shall retain this convention (without explicitly distinguishing the neighborhood and its spine). Note that a grope embedded in a $3$-manifold is automatically framed, so we can be sloppy with the capitalization.
Finally, if one removes a disk from the bottom surface of a grope, one obtains an {\em annular} grope that has two boundary circles. 

\begin{defn} \label{def:grope concordance}
Two links in $S^3$ are called {\em height $n$-Grope concordant} if they cobound an annular height $n$-Grope in $S^3 \times I$.  If an link is height $n$~Grope concordant to the unlink then we say that it is {\em height $n$-Grope slice}. If the links lie in $S^3 \smallsetminus L$ for some other link $L$, and also the annular Grope lies in $(S^3 \smallsetminus L) \times I$, then we refer to it as a Grope concordance {\em rel L}. 
\end{defn}

It is easy to see that the quotient groups $\SC/\SG_n$ from the introduction is the same as the quotient group given by knots, modulo the equivalence relation of height $n$-Grope concordance.

\begin{figure}[htbp]
\setlength{\unitlength}{1pt}
\begin{picture}(216,142)
\put(2,75){$H_1$}
\put(223,75){$H_0$}
\put(195,66){$\theta$}
\put(20,20){\includegraphics{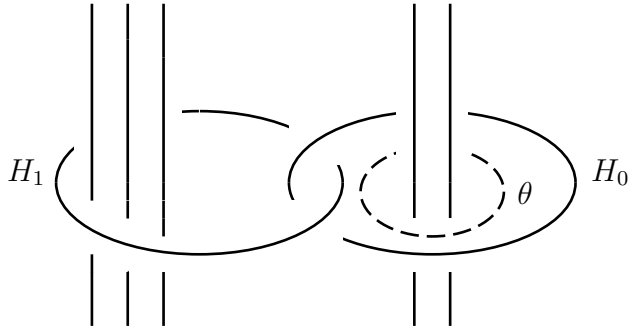}}
\end{picture}
\caption{A Hopf link in $S^3-R$}\label{Hopfsurgery}
\end{figure}

Suppose $R$ is a knot in $S^3$ and $H=\{H_0,H_1\}$ is a $0$-framed
Hopf link in $S^3$ that misses $R$, such as is shown in
Figure~\ref{Hopfsurgery}. Since $0$-framed surgery on a Hopf link
in $S^3$, denoted $(S^3)_H$, is well known to be homeomorphic to
$S^3$ (see \cite{K}), the image of the knot $R$ under this homeomorphism is a new
knot $K$ (i.e. $(S^3,R)_H\approxeq (S^3,K))$. We say that $K$ is
\emph{the result of surgery on $H$}
(see ~\cite{K}). Let $\theta$ denote
the dashed circle shown in Figure~\ref{Hopfsurgery} that is a
parallel of $H_0$ but does not link $H_1$.

\begin{figure}[H]
\begin{center}
\includegraphics{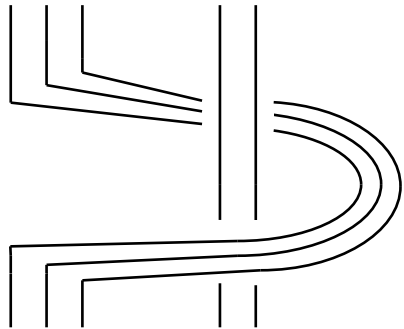}
\caption{The effect of a Hopf surgery}\label{effectofHopf}
\end{center}
\end{figure}

By sliding the left hand strands of $R$ over $H_0$ and then cancelling the Hopf pair, one sees that the effect of one Hopf surgery is as in Figure~\ref{effectofHopf}.

The following proposition discusses the effect of repeated application of the above Hopf surgery. In particular, we assume that several {\em separated} Hopf links are given, i.e. the Hopf links together with their bounding 2-disks are embedded disjointly.

\begin{prop}\label{HopfGrope} Suppose $R$ is a knot and $K$
is the result of surgery on separated Hopf links $H^i$, $1\leq i\leq m$, in $S^3\smallsetminus R$. Suppose the link $\theta_1, \dots, \theta_m$ (see Figure~\ref{Hopfsurgery}) is height $n$-Grope slice rel~$R$. Then $K$ and $R$ are height $n$-Grope concordant.
\end{prop}

\begin{proof}
We first assume that there is only one Hopf link $H=\{H_0,H_1\}$ and later indicate the
modifications necessary in the general case. First observe that, in $S^3$, $K$ and $R$ cobound an embedded punctured annulus. If there is only one strand of $R$ going through $H_1$ in Figure~\ref{Hopfsurgery}, this annulus has one puncture, namely $\theta$, and is shown in Figure~\ref{newannulus}. If there are $r$ strands, just take $r$ parallel copies of this figure, noticing that the relevant part is planar. Thus the result is embedded annulus with $r$ punctures, each a $0$-framed copy of $\theta$. 
\begin{figure}[htbp]
\setlength{\unitlength}{1pt}
\begin{picture}(278,180)
\put(31,115){$R$}
\put(95,141){$K$}
\put(169,90){$\theta$}
\put(-70,-20){\includegraphics{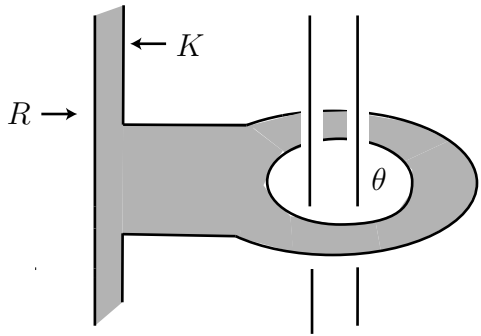}}
\end{picture}
\caption{A punctured annulus cobounding $K$ and $R$}\label{newannulus}
\end{figure}

We can extend this annulus by $R \times  I$ to get a punctured annulus $A$ in $S^3 \times I$, leading from $K$ down to $R$.  
By assumption, there is a height $n$~Grope $G$ that bounds $\theta$ in $(S^3 \smallsetminus R) \times I$. Since $G$ is framed we can get $r$ disjoint parallel copies and we can glue them into the punctures of $A$. The result is a height $n$-Grope concordance between $K$ and $R$. 

In the general case that there is more than one Hopf link, the same arguments works: first one constructs a punctured annulus $A$ between $R$ and $K$, embedded in $S^3$ except for the portion $R \times I$. The punctures of $A$ are now parallels of the link $\theta_1,\dots, \theta_m$ but by assumption they can be filled by disjoint Gropes that miss the interior of $A$ because they lie in the complement of $R \times I$.
\end{proof}

\begin{figure}[htbp]
\setlength{\unitlength}{1pt}
\begin{picture}(213,110)
\put(5,87){$U$}
\put(199,69){$J$}
\put(123,60){$\cong$}
\put(20,10){\includegraphics{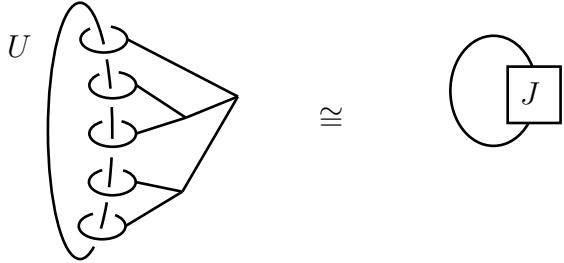}}
\end{picture}
\caption{A clasper construction of $J$}\label{JClasper}
\end{figure}

Now we want to describe a particular large class of knots that
bound $(n+2)$-Gropes in $D^4$ by virtue of satisfying
Proposition~\ref{HopfGrope}. We fix a very specific knot $J$, as
shown implicitly in Figure~\ref{knotJ}, that is obtained
from the trivial knot $U$ by performing $0$-framed surgery on the
$14$-component link shown in
Figure~\ref{knotJ} in the introduction.
For those familiar with the language of
claspers \cite{Hab}, $J$ can also be described as the knot obtained by performing clasper surgery on the unknot along the  height $2$ clasper as shown in
Figure~\ref{JClasper} (we have used a convention wherein an edge corresponds to a left-handed Hopf clasp). It was shown in  \cite{CT}  that $J$ cobounds with the unknot a grope of height~2, embedded in $S^3$. Here we do not want to use this 3-dimensional notion and prefer to give a direct construction of a certain height $(n+2)$-Grope concordance. Note however, that the shift by~2 is a direct consequence of this special feature of $J$.

Suppose $R$ is a slice knot. Given any $\eta_i$, $1\leq i\leq m$,
disjointly embedded circles forming a trivial link in $S^3$ such
that $[\eta_i]\in\1(S^3\smallsetminus R)^{(n)}$, we can consider the knot
$R(\eta_1,\dots,\eta_m,J)$ obtained from $R$ by genetic
infection along each $\eta_i$ by the auxiliary knot $J$ (see
Figure~\ref{infection}).

In the next section we will show that, for many R there exist
certain choices of the \emph{homotopy} classes $[\eta_i]$, such
that for any sufficiently large positive integer $N$, the infected
knot $R(\eta_1,\dots,\eta_m,\#^N_{j=1}J)$ \emph{does not} bound
any Grope of height $(n+2.5)$ in $D^4$. When we say this
we mean for \emph{any} choice of representatives $\eta_i$ of
$[\eta_i ]$. In \emph{this} section, we would like to prove that
$R(\eta_1,\dots,\eta_m,\#^N_{j=1}J)$ \emph{does} bound a 
Grope of height $(n+2)$ in $D^4$. However we cannot prove this in
all cases. We need the (\emph{isotopy} classes) $\eta_i$ to be
chosen carefully to have the stronger {\it geometric} property
that they are height  $n$-Grope slice, rel~$R$, see Definition~\ref{def:grope concordance}. We will demonstrate that this can be arranged by choosing the isotopy
class of each $\eta_i$ carefully within its homotopy class. We
will then use
Proposition~\ref{HopfGrope}, resulting in the following theorem.

\begin{thm}\label{thm:examples} 
Suppose $R$ is a slice knot, $J$ is the knot of Figure~\ref{knotJ}, $N$ is a positive integer,
and $[\eta_i]\in\1(S^3\smallsetminus R)^{(n)}$, $1\le i\le m$.
There exist representatives $\eta_1,\dots,\eta_m$ disjointly
embedded, forming a trivial link in $S^3$, such that
$R(\eta_1,\dots,\eta_m,\#^N_{j=1}J)$ bounds a Grope of
height $(n+2)$ in $D^4$, i.e. is height $(n+2)$-Grope slice.
\end{thm}

\begin{proof} We need the very general Lemma below. Here a $\emph{capped Grope}$ is a Grope equipped with a set of 2-disks (called caps) whose boundaries are
the (full set of) tips of the Grope. In the lemma below, our
capped Gropes are embedded in $S^3$ (although the caps may
intersect the knot $R$). Our only use for the caps is that they
are a good way to make sure that the collection of boundary circles $\eta_i$ of a
disjoint union of such capped Gropes is a trivial link in $S^3$. This follows 
since "one-half" of the caps can be used to ambiently "surger"
the Gropes, producing disjointly embedded disks with boundary
$\eta_i$, each of which lies in a small regular neighborhood of
its corresponding capped Grope. The technique for the following result was used in \cite[Part II Lemma 2.8]{FT}.

\begin{lem}\label{homotopy} Suppose R is a knot and $[\eta_i]\in\1(S^3\smallsetminus R)^{(n)}$, $1\le
i\le m$. Then there exist height $n$  capped Gropes
$G_i$, disjointly embedded in $S^3$, disjoint from $R$ except for
the caps, and such that, for each $i$, $\p G_i$ is in the
homotopy class of $[\eta_i]$.
\end{lem}
\begin{proof}[Proof of Lemma~\ref{homotopy}] The proof is by
induction on $n$. Suppose $n=0$. By general position we can choose
disjoint embedded representatives $\eta_i$ of the $[\eta_i]$.
Moreover, by further ``crossing changes'' we may suppose that
$\eta_i$ forms a trivial link in $S^3$. Setting $G_i=\eta_i$
completes the case $n=0$ since a height zero Grope is merely a
circle. The caps are the set of disks that the $\eta_i$ bound. Now, suppose the theorem is true for $n-1$. Suppose
$[\eta_i]=\prod_j[a_{ij},b_{ij}]$ for some $a_{ij}$, $b_{ij}$ in
$\1(\sk)^{(n-1)}$. Using the induction hypothesis applied to each of the $a_{ij}$ and $b_{ij}$, choose
 embedded height $(n-1)$ capped Gropes $H_{ij}$ and
$L_{ij}$, pairwise disjoint except for the basepoint, such that
$[\p H_{ij}]=a_{ij}$, $[\p L_{ij}]=b_{ij}$ for all $i$ and all
$j$. Let $A_{ij}=\p H_{ij}$, $B_{ij}=\p L_{ij}$. We can alter
these Gropes to assume that, for each fixed $i$, the $A_{i1}$,
$B_{i1}$, $A_{i2}$, $B_{i2},\dots$ share a common point as illustrated for the case $i,j\leq 2$ on the left-hand side of Figure~\ref{wedge}, with the $H_{ij}$ coming up straight out of the
plane of this page and the $L_{ij}$ going down straight below the
plane of this page.
\begin{figure}[htbp]
\setlength{\unitlength}{1pt}
\begin{picture}(178,115)
\put(-8,71){$B_{11}$}
\put(-9,51){$B_{21}$}
\put(5,92){$A_{11}$}
\put(16,24){$A_{21}$}
\put(48,23){$A_{22}$}
\put(47,93){$A_{12}$}
\put(66,44){$B_{22}$}
\put(64,74){$B_{12}$}
\put(10,10){\includegraphics{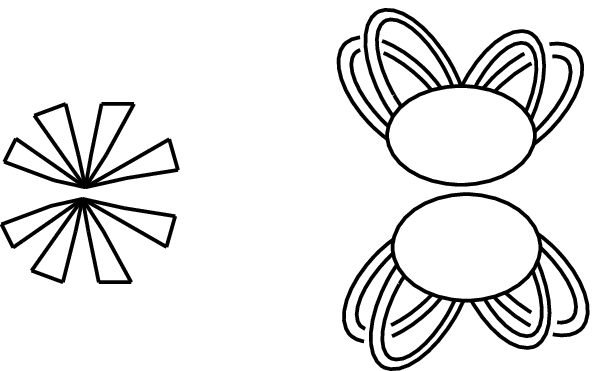}}
\put(134,74){$N_1$}
\put(136,42){$N_2$}
\end{picture}
\caption{}\label{wedge}
\end{figure}

Then we can ``thicken'' this wedge slightly as
shown on the right side of Figure~\ref{wedge}, avoiding all the caps
for $H_{ij}$ and $L_{ij}$, to form a new embedded surface $N_i$
whose boundary has the homotopy class of $[\eta_i]$. Then
$N_i\cup(\bigcup_jH_{ij}\bigcup_jL_{ij})$ forms a height $n$ capped
Grope $G_i$, and these are disjoint for different $i$.
\end{proof}

\begin{rem} In the above theorem, $S^3$ may be replaced by any
orientable $3$-manifold $X$ as long as the $\eta_i$ are trivial in
$\1(X)$.
\end{rem}

Continuing with the proof of Theorem~\ref{thm:examples}, by
Lemma~\ref{homotopy} there exist embedded representatives
$\eta_i$, $1\le i\le m$, forming a trivial link in $S^3$, that
bound disjointly embedded height $n$ Gropes in $S^3\smallsetminus R$
(with caps that intersect $R$). We get height $n$-Gropes $G_1,\dots,G_m$ in $(S^3 \smallsetminus R) \times I$ by pushing slightly into the $I$-direction, leaving only the boundary in $S^3$. 

For each $i$, form $5N$ parallel push-offs of $\eta_i$, denoted $\eta_{ijk}$, $1\le j\le N$, $1\le
k\le 5$. Now, for each fixed $i$ and $j$, connect the $5$ circles by four short arcs, connecting
$\eta_{ij1}$ to $\eta_{ij2}$,  $\eta_{ij2}$ to $\eta_{ij3}$, et cetera. For each fixed $i$ and $j$, \emph{inside a regular neighborhood of the short arcs}, connect the $5$ circles
$\eta_{ij1},\dots,\eta_{ij5}$ by a collection of arcs to form a copy of the
embedded tree as shown in the left hand side of Figure~\ref{JClasper}. Then for varying $i,j$, using the $mN$ trees as guides, replace each tree with a copy, $L_{ij}$, of the
$14$-component link of Figure~\ref{knotJ}. These
trees, and hence the links, will be pairwise disjoint. Let $K$ denote the result of $0$-framed surgery on $\coprod_{i,j}L_{ij}$.

First we claim that $K$ is isotopic to
$R(\eta_1,\dots,\eta_m,\#^N_{j=1}J)$. For, since for fixed
$(i,j)$, $\eta_{ij1},\dots,\eta_{ij5}$ are parallel, there is a
3-ball $B_{ij}$ in $S^3$ that contains $L_{ij}$  and such that
the pair $(B_{ij}, B_{ij} \bigcap R)$ is a trivial tangle. In this
way one sees that the effect of the $ij^{th}$ surgery is locally
the same as the effect of the $14$-component link of
Figure~\ref{knotJ} on a trivial tangle. This effect is clearly to
tie all the parallel strands of the trivial tangle into the knot
$J$. Moreover, since for each fixed $i$, the $N$ circles
$\eta_{ij1}$, $1\le j\le N$ are also parallel, modifying each by
$J$ is the same as infecting a single one, say $\eta_{i11}$, by
$\#^N_{j=1}J$. Thus $K$ is the result of genetic modification of
$R$ along the original $m$ circles $\{\eta_{i}\}$ via
$\#^N_{j=1}J$.

Secondly we claim that $K$ satisfies the hypotheses of Proposition
~\ref{HopfGrope} for $(n+2)$ and so $K$ and $R$ are height $(n+2)$-Grope concordant. Since $R$ was assumed slice, the verification of this claim
will finish the proof of the theorem. To apply that Proposition we
must establish that $K$ can be viewed as the result of surgery on
a disjoint union of Hopf links in the exterior of some slice knot
and show that the appropriate circles $\theta$ associated to each Hopf link are height $(n+2)$-Grope slice rel~$R$. We focus on one link $L_{ij}$ and order its components $(H_0,\dots,H_{13})$ as shown in Figure~\ref{JGrope}. Ignore the shading for now. 
\begin{figure}[htbp]
\setlength{\unitlength}{1pt}
\begin{picture}(242,190)
\put(198,150){$H_0$}
\put(50,138){$H_1$}
\put(231,70){$H_2$}
\put(190,106){$H_3$}
\put(201,34){$H_4$}
\put(133,13){$H_5$}
\put(118,38){$H_6$}
\put(129,66){$H_7$}
\put(127,111){$H_8$}
\put(178,75){$H_9$}
\put(49,25){$H_{10}$}
\put(50,55){$H_{11}$}
\put(50,78){$H_{12}$}
\put(50,110){$H_{13}$}
\put(10,10){\includegraphics{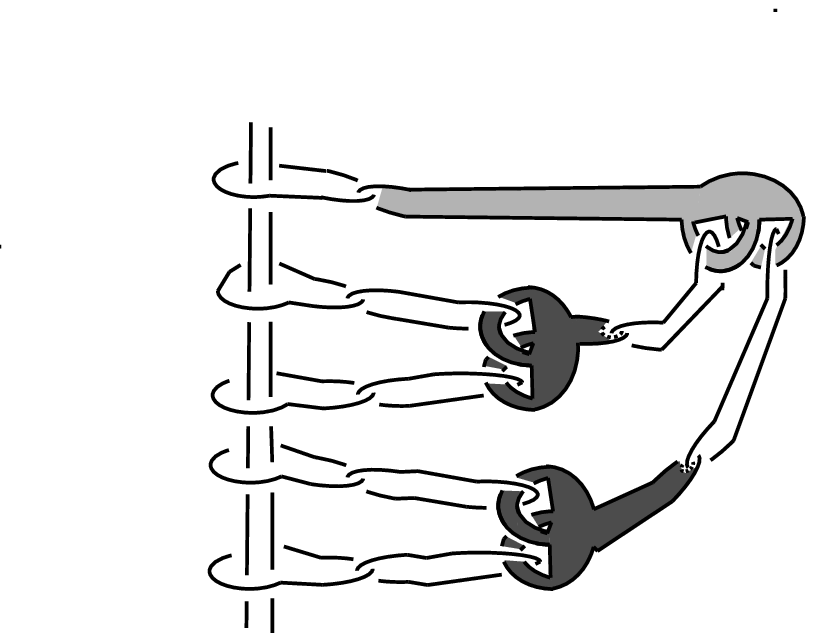}}
\end{picture}
\caption{The link $L_{ij}$}\label{JGrope}
\end{figure}

The pictured link is merely a copy of that shown in Figure~\ref{knotJ}, except that the three copies of the Borromean rings have been altered by isotopy to arrive at the embedding shown in Figure~\ref{JGrope}. We can now perform 4 handle cancellations, namely the 4 pairs on the left, read from the bottom as
\[
H_{10}-H_5, H_{11}-H_6, H_{12}-H_7, H_{13}-H_8.
\]
What remains is a 6 component link $N_{ij}$ shown in Figure~\ref{revisedlinkLij}.
\begin{figure}[htbp]
\setlength{\unitlength}{1pt}
\begin{picture}(242,190)
\put(193,145){$H_0$}
\put(55,137){$H_1$}
\put(231,71){$H_2$}
\put(192,108){$H_3$}
\put(196,33){$H_4$}
\put(167,83){$H_9$}
\put(10,10){\includegraphics{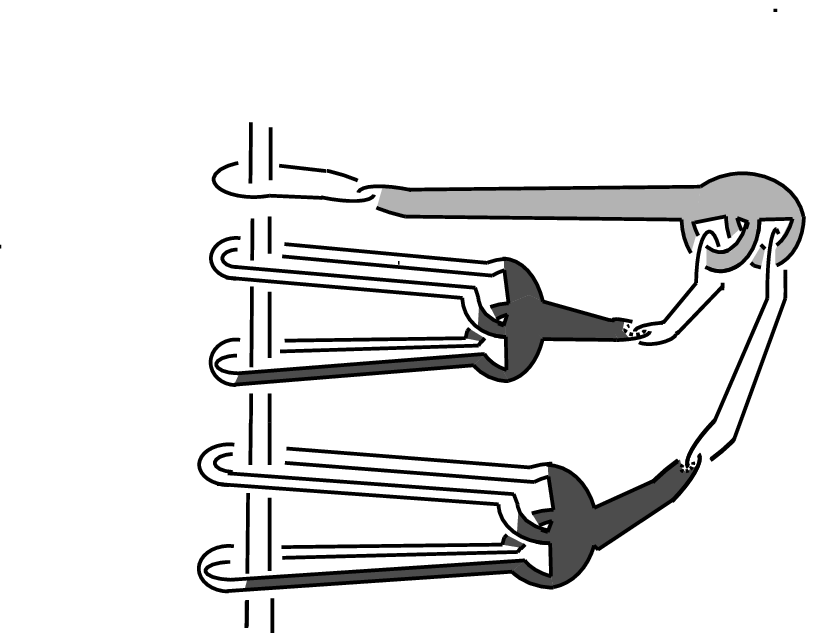}}
\end{picture}
\caption{The link $N_{ij}$}\label{revisedlinkLij}
\end{figure}

 Let $L:=N_{ij}-\{H_0,H_1\}$.  Since $L$ consists of two separated Hopf links, surgery along $L$ yields $S^3$ again. Furthermore, the knot $R$ can be disentangled from these two Hopf links by sliding over $H_2$ and $H_3$. Therefore, the image of $R$ in $S^3_L$ remains a slice knot, $R'$. Thus we have shown that $K$ is the image of a slice knot,$R'$, in a manifold, $S^3_L$, homeomorphic to $S^3$, after surgery on a $2$-component link $(H_0,H_1)$. Using handle slides over $H_4$ and $H_9$ one sees that this link is indeed a Hopf link in $S^3_L$. 
 
We now describe a height $2$ Grope, embedded in $S^3$ whose boundary is $\theta$. Recall from the notation of Figure~\ref{Hopfsurgery} that $\theta$ is merely a parallel of $H_0$ that is short-circuited to avoid linking $H_1$. The (punctured) torus shown shaded lightly in
Figure~\ref{revisedlinkLij} is the first stage of the Grope. The obvious symplectic basis for this torus consists of two circles that are (isotopic to) meridians of $H_2$ and
$H_3$, shown as dashed in Figure~\ref{revisedlinkLij}. The annuli that go between the symplectic basis and the dashed circles are part of the Grope and are often called ``pushing annuli''. These dashed meridians in turn bound punctured torii, that are shown darkly shaded (partially) in Figure~\ref{revisedlinkLij}.
This gives the embedded Grope $G$ of height $2$,
each of whose surfaces has genus one, together with four pushing annuli identifying the tips
of the Grope with $\{\eta_{2},\eta_{3},\eta_{4},\eta_{5}\}$. Recall that we are fixing $i,j$ and hence there are really $mN$ such Gropes $G_{ij}$ with boundary $\theta_{ij}$, embedded disjointly in $S^3_L$. 

The Gropes $G_{ij}$ will constitute the first two stages of the height $(n+2)$ Gropes whose boundary is $\theta_{ij}$, that we need to exhibit in $(S^3_L\smallsetminus R) \times I$ to verify the hypothesis of Proposition~\ref{HopfGrope}. The higher stages of this Grope are just $4N$ parallel copies of the height $n$-Gropes $G_i$ bounding $\eta_i$ that exist by assumption. These are attached to the tips of $G_{ij}$.  This finishes the verification that $K$ satisfies the hypotheses of Proposition~\ref{HopfGrope} for $(n+2)$ and so $K$ bounds a Grope in
$D^4$ of height equal to $(n+2)$.
\end{proof}

We note in passing that the argument used in the verification that
the knot $K$ satisfied the hypotheses of Proposition
~\ref{HopfGrope}, generalizes to establish the following Corollary
(that we make no use of in this paper). We state it in the
language of claspers and the trees that underly them. See [CT] for
definitions. It has generalizations to non-symmetric trees but we
suppress this for simplicity. Suppose $(T,r)$ is a uni-trivalent
tree $T$ equipped with a distinguished univalent vertex $r$ called
the {\it root}. Each (tree) clasper has an underlying tree (compare Figure~\ref{JClasper}).

\begin{cor}\label{n-Grope} Suppose $K$ is
the result of clasper surgery on a knot $R$ along
$\coprod^\ell_{i=1}C_{(T_i,r_i)}$ (a disjoint union of claspers
whose root leaves bound disjoint disks in $S^3$). Suppose each
$(T_i,r_i)$ is a rooted symmetric tree of height $m$ and that the 
non-root leaves of $C_{(T_i,r_i)}$ form a link that is height $n$-Grope slice rel~$R$. Then $K$
and $R$ are height $(m+n)$-Grope concordant.
\end{cor}

\section{Knots that are not \nss}\label{sec:not solvable}

In this section we complete the proof of our main Theorem~\ref{thm:main}. First we review the necessary definitions.

For a CW-complex $W$, we define $W^{n}$ to be the regular
covering corresponding to the subgroup $\pi_1(W)^{(n)}$. If $W$ is an
oriented $4$-manifold then there is an intersection form
$$
\lambda_{n}:H_2(W^{(n)}) \times H_2(W^{(n)})
\lra\bbz[\pi_1(W)/\pi_1(W)^{(n)}].
$$
Similarly there is a self-intersection form $\mu_n$. For $n\in\bbn_0$, an
{\it $n$-Lagrangian} is a submodule $L\subset H_2(W^{(n)})$ on which
$\lambda_n$ and $\mu _n$ vanish and which maps onto a Lagrangian of
$\lambda_0$.

\begin{defn} \cite[Sections 7,8]{COT} A knot is called {\it $(n)$-solvable} if $M$ (the zero
surgery on the knot $K$) bounds a spin $4$-manifold $W$, such
that the inclusion map induces an isomorphism on first
homology (a $4$-manifold satisfying only these is called an \emph{$H_1$-bordism} for $M$) and such that W admits two {\it dual}
$n$-Lagrangians. This means that the intersection form $\lambda_n$ pairs
the two Lagrangians non-singularly and that their images
together freely generate $H_2(W)$. Then $M$ is also called $(n)$-solvable and $W$ is called an \emph{$(n)$-solution} for $M$ and $K$. A knot is called {\it $(n.5)$-solvable}, $n\in\bbn_0$, if $M$ bounds a
spin $4$-manifold $W$ such that the inclusion map induces an isomorphism
on first homology and such that $W$ admits an $(n+1)$-Lagrangian and a
dual $n$-Lagrangian in the above sense. The $W$ is called an \emph{$(n.5)$-solution} for $K$ and $M$.
\end{defn}

Recall also the \emph{rational derived series}, $G^{(n)}_r$, of a group $G$ wherein $G^{(0)}_r\equiv G$ and $G^{(n+1)}_r\equiv \{g \mid g^k\in [G^{(n)}_r,G^{(n)}_r]$ for some positive integer $k\}$. The terms of the rational derived series are slightly larger than the terms of the derived series but they have the key technical advantage that their successive quotients are torsion-free abelian groups \cite[Section 3]{Ha}.

 The following are the major theorems of this section. Combined with Theorem~\ref{thm:examples} they give large families of knots that are $(n)$-solvable but not
 $(n.5)$-solvable. 

\begin{thm}\label{thm:infection} Let $R$ be an \ns\ knot
($n\ge1$) and $M$ the $0$-framed surgery on $R$. Suppose there
exists a collection of homotopy classes 
\[
[\eta_i]\in\1(M)\2n, \quad 1\leq i\leq m, 
\]
that has the following property: For \emph{any}
$(n)$-solution $W$ of $M$ there exists {\bf some} $i$ such that
$j_*(\eta_i)\notin\1(W)\nn_r$ where $j_*:\1(M)\lra\1(W)$. 

Then, for any oriented trivial link $\{\eta_1,\dots,\eta_m\}$ in $S^3\smallsetminus R$
that represents the $[\eta_i]$, and for any $m$-tuple
$\{J_1,\dots,J_m\}$ of Arf invariant zero knots for which
$|\rho_\BZ(J_i)|>C_M$ (the Cheeger-Gromov constant of M), the knot
\[
K=R(\eta_1,...,\eta_m,J_1,...,J_m)
\]
formed by genetic modification is \ns\ but not \nss. Moreover, $K$ is of infinite
order in $\ff$.
\end{thm}

\begin{thm}\label{thm:injectivity} Suppose $R$ is a genus 2
fibered knot that is \ns\ (for example a genus 2 fibered ribbon
knot). Then there exists a collection of homotopy classes
satisfying the hypotheses of Theorem~\ref{thm:infection}.
\end{thm}

\begin{proof}[Proof of Theorem~\ref{thm:main}] First, we will show that $\ff$ contains an element of infinite order.
 We may assume that $n\geq 1$ since this result was known previously for $n=0,1$ and $2$. There exist genus 2 fibered ribbon knots, for example the connected sum
of two figure eight knots. Hence by Theorem~\ref{thm:injectivity}
there exists an \ns\ knot $R$ and classes $\{[\eta_i]\}$
satisfying the hypotheses of Theorem~\ref{thm:infection}. Certainly
there exist representatives of these classes that form a trivial
link since we can alter any collection by crossing changes to
achieve this. There also exist Arf invariant zero knots with
$\rho_\BZ>C_M$, for example the connected sum of a suitably large
{\bf even} number of left-handed trefoil knots. Then
Theorem~\ref{thm:infection} implies that the knot $K$ formed by
genetic modification of $R$ along any such collection $\{\eta_i\}$
using any such $\{J_i\}$ is of infinite order in $\ff$.

Now, in order to prove the other statements of Theorem~\ref{thm:main}, we just need to be a little more careful in choosing the infection knots. In fact, the shift by two for the inclusion $\SG_{n+2}\leq \SF_n$ is related to the fact that the knot $J$ from Figure~\ref{knotJ} by construction cobounds with the unknot a Grope of height~2 in $S^3$, yet it also satisfies the following essential nontriviality condition. (All that we truly require is that the integral, over the circle, of the Levine-Tristram signature function of $J$ is \emph{non-zero}.) 

\begin{lem}\label{signatureJ}
The von Neumann signature of $J$ associated to the infinite cyclic cover, denoted $\rho_{\mathbb{Z}}(J)$, satisfies
\[
\rho_{\bbz}(J) =\frac {4}3
\]
In fact, $J$ has the same Levine-Tristram signatures as the left-handed trefoil knot.
\end{lem}

We shall prove this Lemma at the end of this section.

\begin{proof}[Proof of Theorem~\ref{thm:infection}] Let $V'$ be
an $(n)$-solution for $M$. Let $N$ be the zero surgery on
$K$. Using $V'$ we will show that $N$ (and hence $K$) is \ns.
Moreover, assuming the existence of an $(n.5)$-solution $V$
for $N$, we shall derive a contradiction.

There is a {\sl standard cobordism} $C$ between $M$ and $N$ which can be described as follows. For each $1\le
i\le m$, choose a spin $4$-manifold $W_i$ whose boundary is
$M_{J_i}$, the zero surgery on $J_i$, such that $\1(W_i)\cong\BZ$
generated by a meridian of $J_i$ and such that the intersection
form on $H_2(W_i)$ is a direct sum of hyperbolic forms. Such a
$W_i$ is a $0$-solution of $M_{J_i}$. It exists whenever the Arf
invariant of $J_i$ is zero. Now form $C$ from $M\x[0,1]$ and
$\coprod^m_{i=1}W_i$ by identifying (for each $i$) the solid torus
in $\p W_i\equiv(S^3\backslash J_i)\cup (S^1\x D^2)$ with a
regular neighborhood of $\eta_i\x\{1\}$ in such a way that a
meridian of $J_i$ is glued to a longitude of $\eta_i$ and a
longitude of $J_i$ is glued to a meridian of $\eta_i$. Then
$\p_+C\cong N$ and $\p_-C\cong M$. Also observe that $C$ can be assumed to be spin by changing the spin structure on $W_i$ if necessary.

Let $W=C\cup V$ and $W'=C\cup V'$. We claim that $W$ is an
$(n)$-solution for $M$ and that $W'$ is an $(n)$-solution for $N$.
Clearly $W$ is an $H_1$-bordism for $M$ and $W'$ is an
$H_1$-bordism for $N$, and 
\[
H_1(M)\cong H_1(C)\cong H_1(N)\cong
H_1(V)\cong H_1(W)\cong H_1(W')\cong\BZ
\]
 all induced by inclusion.
Thus $V$ and $N=\p V$ have 2 distinct spin structures and changing
the spin structure on $V$ changes that induced on $\p V$
(similarly for $M$ and $V')$. Hence spin structures on the
manifolds $V$ and $V'$ can be chosen to agree with those induced
on $N$ and $M$ by the spin manifold $C$. Thus $W$ and $W'$ are
spin. It is easy to see that $H_2(C)$ is isomorphic to
$\bigoplus^m_{i=1}H_2(W_i)\op H_2(M)$ by examining the
Mayer-Vietoris sequence for
$C\cong(M\x[0,1])\cup(\prod^m_{i=1}W_i)$. By a similar sequence
for $W\cong V\cup C$ one sees that 
\[
H_2(W)\cong (H_2(C)\op
H_2(V))/i_*(H_2(N)).
\]
Since $N\lra V$ induces an isomorphism on
$H_1$, one easily sees that it induces the zero map on $H_2$.
Moreover, a generator of $H_2(N)$ under the map 
\[
H_2(N)\lra
H_2(C)\cong\bigoplus^m_{i=1}H_2(W_i)\op H_2(M)
\]
 goes to a
generator of $H_2(M)$ since it is represented by a capped-off
Seifert surface for $R$ that can be chosen to miss the $\eta_i$
(since $n\ge1$). It follows that 
\[
H_2(W)\cong
\bigoplus^m_{i=1}H_2(W_i)\op H_2(V).
\]
 Similarly, $H_2(W')\cong
\bigoplus^m_{i=1}H_2(W_i)\op H_2(V')$. Since $V$ is an
$(n.5)$-solution for $N$, it is an $(n)$-solution for $N$, so there
exists an $n$-Lagrangian with $n$-duals. This may be thought of as
collections $\SL$ and $\SD$ of based immersed surfaces that lift
to the $\1(V)\2n$ cover of $V$ (called $n$-surfaces in
\cite[Sections 7 and 8]{COT}), that together constitute a basis of
$H_2(V;\BZ)$ and that satisfy
\[
\la_n(\ell_i,\ell_j)=\mu_n(\ell_i)=0, \quad \la_n(\ell_i,d_j)=\d_{ij}
\]
 for the intersection and
self-intersection forms with $\BZ[\1(V)/\1(V)\2n]$-coefficients.
These same collections are certainly $n$-surfaces in $W$ since
$\1(V)\2n$ maps into $\1(W)\2n$ under the inclusion. Since the
equivariant intersection form can be calculated from the
collections $\SL$ and $\SD$, (or by naturality) these also retain
the above intersection properties with $\BZ[\1(W)/\1(W)\2n]$
coefficients. Now consider collections of $0$-Lagrangians, $\SL_i$, and $0$-duals, $\SD_i$, for the
$0$-solutions $W_i$. Since the map
\[
\1(W_i)\lra\1(C)\lra\1(W)/\1(W)\2n
\]
 is zero (since $\1(W_i)$ is
generated by $\eta_i$), these $0$-surfaces are $n$-surfaces in
$W$. Since the union of all of these collections is a basis for
$H_2(W)$, it constitutes an $n$-Lagrangian and $n$-duals for $W$.
Thus $W$ is an $(n)$-solution for $M$. Similarly, $W'$ is an
$(n)$-solution for $N$. In particular $K$ is \ns.

Now let $\G=\1(W)/\1(W)\nn_r$. It is straightforward to verify
that $\G$ is an  \ns\ poly-(torsion-free-abelian) group (abbreviated PTFA) \cite[Section 3]{Ha}. Let
$\psi:\1(W)\lra\G$ be the projection. Let $\phi$ and $\phi'$
denote the induced maps on $\1(M)$ and $\1(N)$ respectively.

\begin{prop}\label{rho invariants}
$$
\rho_\G(M,\phi) - \rho_\G(N,\phi') =
\sum^m_{i=1}\e_i\rho_\BZ(J_i)
$$
where $\e_i=0$ or $1$ according as $\phi(\eta_i)=e$ or not.
\end{prop}

\begin{proof} By \cite[Proposition 3.2]{COT2}
$$
\rho_\G(M,\phi) - \rho_\G(N,\phi') = \sum^m_{i=1}
\rho_\G\(M_{J_i},\psi_{|\1(M_{J_i})}\).
$$
Since $\psi_{|\1(M_{J_i})}$ factors through $\BZ$ generated
by $\eta_i$, its image is zero if $\e_i=0$ and is $\BZ$ if
$\e_i=1$. In the former case $\rho_\G=0$ since $\s^{(2)}_\G$
is then the ordinary signature. In the latter case
$\rho_\G=\rho_\BZ$ (replacing $\G$ by the image of
$\psi_{|\1(M_{J_i})}$), by
\cite[Proposition 5.13]{COT}. But
$\rho_\BZ(M_{J_i},\psi_{|\1(M_{J_i})})$ is just what we have
called $\rho_\BZ(J_i)$.
\end{proof}

Moreover since $V$ is assumed to be an
$(n.5)$-solution for $N$ and $\phi'$ extends to $\1(V)$,
$\rho_\G(N,\phi')=0$ by \cite[Theorem 4.2]{COT}. Finally by
hypothesis, there exists some $i$ such that
$\phi(\eta_i)\neq e$. Thus, by Proposition~\ref{rho invariants}, $|\rho_\G(M,\phi)|>C_M$, a
contradiction. Therefore no $(n.5)$-solution exists for $N$, that is to say $K$ is not $(n.5)$-solvable.

Now suppose $K=K(R,\eta,J)$ were of finite order $k>0$ in \ff. Let
$N_\#$ denote the $0$-framed surgery on $\#^k_{j=1}K$, and
let $V_\#$ denote an $(n.5)$-solution for $N_\#$. We shall arrive
at a contradiction, implying that $K$ is of infinite order.
There is a standard cobordism $D$ from $\coprod^k_{j=1}N_j$
(where $N_j$ is the $j^{\supth}$ copy of $N$) to $N_\#$ that
is obtained from $\coprod^k_{j=1}N_j\x[0,1]$ by adding $(k-1)$
1-handles and then $(k-1)$ 2-handles. Let $W$ denote the
4-manifold 
\[
W:=V\cup D\cup\(\coprod^k_{j=1}C_j\)
\]
where $C_j$ is
the $j^{\supth}$ copy of the standard cobordism $C$
(described above) between $M$ and $N$. Let
$\G=\1(W)/\1(W)^{(n+1)}_r$ and $\psi:\1(W)\to\G$ the
projection. Let $\phi_j$, $\phi'_j$, $\phi$ denote the
restrictions of $\psi$ to $\1(M_j)$, $\1(N_j)$ and $\1(N_\#)$
respectively. Since $D$ is a homology cobordism with
coefficients in $\SK(\G)$, the quotient field of $\BZ\G$
\cite[Lemma 4.2]{COT2}, $\sum^k_{j=1}\rho(N_j,\phi'_j)=
\rho(N_\#,\phi_\#)$. Since $V_\#$ is an $(n.5)$-solution for
$N_\#$ over which $\phi_\#$ extends, $\rho(N_\#,\phi_\#)=0$
by \cite[Thm. 4.2]{COT}. Thus,
\[
\sum^k_{j=1}\rho(M_j,\phi_j)=\sum^k_{j=1}\sum^m_{i=1}\e_{ji}\rho_\BZ(J_i)
\]
where $\e_{ji}=0$ if $\phi_j(\eta_{ji})=e$ and $\e_{ji}=1$ if
$\phi_j(\eta_{ji})\neq e$. We assert (and prove below) that,
for each $j$, there exists some $i$, such that
$\phi_j(\eta_{ji})\neq e$. Thus
$|\sum^k_{j=1}\rho(M_j,\phi_j)|>kC_M$, which is a contradiction
since $|\rho(M,\phi_j)|<C_M$. To establish the assertion
above, fix $j=j_0$ and consider the $4$-manifold 
\[
\ov W=W\cup\(\bigcup^k_{\substack{j=1\\ j\neq j_0}}-W_j\)
\]
 where
$W_j$ is a copy of the $(n)$-solution for $M$ discussed above.
Then $\p\ov W=M_{j_0}$. By an argument as above, one sees that $\ov W$ is an $(n)$-solution for $M_{j_0}$. By the assumption in Theorem~\ref{thm:infection}, there
exists some $i$ (depending on $j_0$) such that
$j_*(\eta_{j_0i})\notin\1(\ov W)^{(n+1)}_r$ where
$j:M_{j_0}\to\ov W$ is the inclusion. Since this inclusion
factors $M_{j_0}\overset{i_*}{\lra}W\lra\ov W$, it follows
that $i_*(\eta_{j_0i})\notin\1(W^{(n+1)}_r)$. Hence
$\phi_{j_0}(\eta_{j_0i})\neq e$ in $\G$. This establishes the
assertion and completes the proof of Theorem~\ref{thm:infection}.
\end{proof}

\begin{proof}[Proof of Lemma~\ref{signatureJ}]
Recall  we can compute $\rho_{\mathbb{Z}}(J)$ from any $\mathbb{Z}$-bordism $W$ for $M_J$ as follows. Let $M$ be a matrix for the intersection form on $H_2(W;\bbz[\bbz]))$. Then the integral over the unit circle of the signature of $M(\omega)$ gives the von Neumann signature of $W$ associated to the map to $\mathbb{Z}$, $\sigma_{\mathbb{Z}}(W)$, and by definition $\rho_{\mathbb{Z}}(J)=\sigma_{\mathbb{Z}}(W)-\sigma (W)$.
Recall that $J$ is obtained from the unknot $U$ by $0$-framed surgery on the $14$-component link shown in Figure~\ref{knotJ} (see also Figure~\ref{JGrope} for the labelling of the components). This picture also encodes an assortment of possible 4-manifolds $W$, each with boundary the $0$-surgery $M_J$. This is seen as follows. Consider labelling some of the components of the $14$-component link by dots and label the remainder by $0$'s. Also label the unknot $U$ by a dot, subject to the restriction that the collection, $A$, of all components labelled with dots constitutes a trivial link in $S^3$, hence bounds a collection of disjointly embedded 2-disks in the $4$-ball. Let $V$ be the exterior of this collection of disks in the $4$-ball. Now interpret the remaining circles with $0$ labels as $0$-framed $2$-handles attached to $S^3-A$ which is part of the boundary of $V$. The result is a $4$-manifold, $W$.

Since $V$ is diffeomorphic to the manifold obtained from a $4$-ball by \emph{adding} $1$-handles (one for each dotted component), one can think of the choices above as deciding which components of the $14$-component link will represent $1$-handles and which should represent $2$-handles. Also since $V$ has the same \emph{boundary} as the manifold obtained by adding $0$-framed 2-handles along the trivial link represented by $A$, the boundary of \emph{any} such $W$ is the same as doing $0$-framed surgery on \emph{all} the components, namely $0$-framed surgery on $J$. 
\begin{figure}[htbp]
\setlength{\unitlength}{1pt}
\begin{picture}(251,183)
\put(37,176){$U$}
\put(198,150){$H_0$}
\put(53,142){$H_1$}
\put(231,70){$H_2$}
\put(192,110){$H_3$}
\put(202,37){$H_4$}
\put(133,17){$H_5$}
\put(118,43){$H_6$}
\put(129,70){$H_7$}
\put(127,118){$H_8$}
\put(178,77){$H_9$}
\put(50,31){$H_{10}$}
\put(50,59){$H_{11}$}
\put(50,81){$H_{12}$}
\put(50,111){$H_{13}$}
\put(10,10){\includegraphics{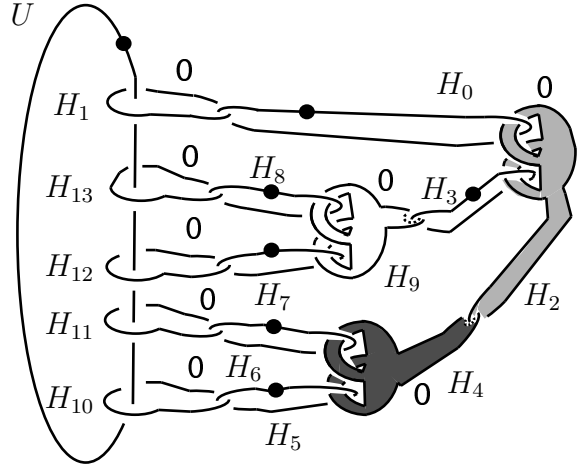}}
\end{picture}
\caption{A 4-manifold with boundary $M_J$}\label{fourmanifoldW}
\end{figure}

The best decision for our purposes is as shown in Figure~\ref{fourmanifoldW}. Note that the Borromean rings formed by $\{H_0,H_2,H_3\}$ has been altered by an isotopy to yield a different embedding than that shown in Figure~\ref{JGrope}. This allows for six $1-2$-handle cancellations, so that the diagram can be simplified by handle slides to leave only one $1$-handle (corresponding to $U$) and two $2$-handles, corresponding to $H_2$ and $H_4$. The resulting handle decomposition is shown in Figure~\ref{fourmanifoldW3}. 
\begin{figure}[htbp]
\setlength{\unitlength}{1pt}
\begin{picture}(248,179)
\put(42,150){$U$}
\put(154,153){$H_2$}
\put(105,17){$H_4$}
\put(200,40){$0$}
\put(232,151){$0$}
\put(10,10){\includegraphics{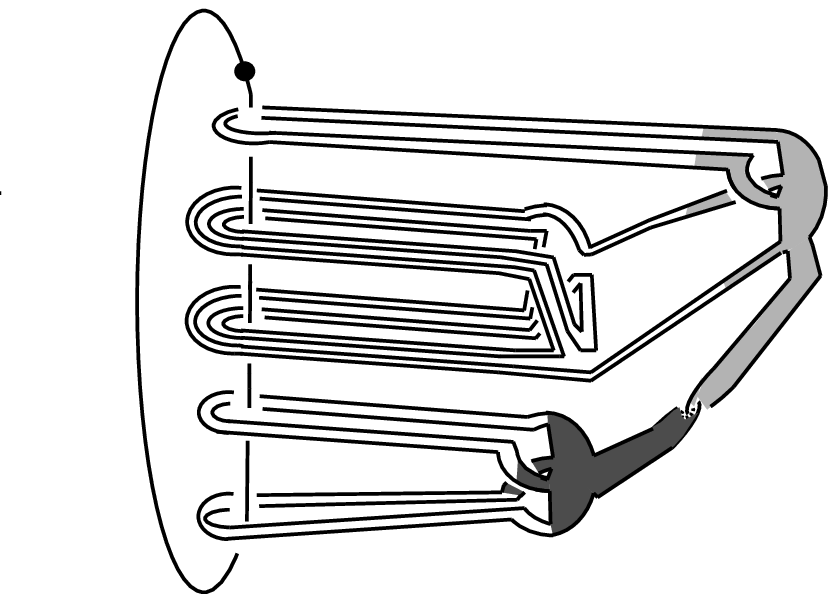}}
\end{picture}
\caption{}\label{fourmanifoldW3}
\end{figure}

It follows that $\1W \cong \bbz$, generated by a meridian, and that $\pi_2W$ is free on the two $2$-handles $H_2,H_4$ in the sense that the core of a 2-handle leads to a $2$-sphere in $W$ if the attaching circle is null homotopic in the 4-ball (minus the disk corresponding to the 1-handle).  In fact, the circle labelled $H_4$ bounds an immersed $2$-disk that can be seen as follows. Note that $H_4$ bounds the(punctured) torus shown darkly shaded Figure~\ref{fourmanifoldW3}
But the two circles representing a symplectic basis for this torus are freely homotopic meridians of $U$. Thus the circle labelled $H_4$ in Figure~\ref{fourmanifoldW3} can be altered by a local homotopy (no intersections with $U$) until it is undone from $U$. Similarly, $H_2$ winds around $U$ in a way that is precisely a ``doubling'' of the way $H_4$ wraps around $U$ in Figure~\ref{fourmanifoldW3}. It follows that $H_2$ is homotopic to a circle that is undone from $U$. This implies
\[
W\simeq S^1 \vee S^2 \vee S^2.
\]
 Examining the intersections created during these homotopies (in the group ring $\bbz[\bbz]$) gives the intersection form $ \lambda$ on $\pi_2W$. For example, since the two homotopies above are disjoint, there are immersed $2$-spheres representing $[H_2]$ and $[H_4]$ that intersect in precisely one point. Therefore $\lambda([H_2],[H_4])=1$. It turns out that the intersection matrix is given by
\[
\lambda = \left( \begin{matrix}
-t-t^{-1}+2 &1  \\
1 & -t-t^{-1}+2 \\
\end{matrix} \right)
\]
Since the determinant $(t+t^{-1}-2)^2 -1$ of $ \lambda$ factors as ($t+t^{-1}-1$)($t-3+t^{-1}$)), its only zeroes on the unit circle are the two primitive $6^{th}$ roots of unity. These are exactly the zeros of the first factor, which is the Alexander polynomial of the trefoil knot. Since the signature function on the unit circle is constant except at these values, it suffices to check that the signature of $ \lambda$ is $+2$ when substituting $t=-1$. Note that if both diagonal entries of $\lambda$ were changed in sign, then the analysis is the same except that $J$ has the signature of the right-handed trefoil knot. Since all we really need is that the integral of the signature function is \emph{non-zero}, this would also suffice.
\end{proof}

We can now finish the proof of Theorem~\ref{thm:main}.
Let $R$ be a genus 2 fibered ribbon knot. Let $M$ be the zero
surgery on $R$ and let $N$ be an even integer greater than $C_M$
where $C_M$ is the \textit{Gromov-Cheeger constant} for $M$. Given
$n$, by Theorem~\ref{thm:injectivity} there exist homotopy classes
$[\eta_i]\in\1(\sr)^{(n)}$, $1\le i\le m$, with the property that,
for any $(n)$-solution $W$ of $M$ there exists some $i$ such that
$[\eta_i]$ is not mapped to zero under the homomorphism
$\1(M)\to\1(W)\to\1(W)/\1(W)^{(n+1)}_r$.  Let $J$ be the knot
shown in Figure~\ref{knotJ}. Applying Theorem~\ref{thm:examples} we
get representatives $\eta_1,\dots,\eta_m$ disjointly embedded,
forming a trivial link in $S^3$ such that
$K=R(\eta_1,\dots,\eta_m,\#^N_{j=1}J)$ bounds a Grope of
height $(n+2)$ in $D^4$. Hence
$R(\eta_1,\dots,\eta_m,\#^N_{j=1}J)\in \SG_{(n+2)}$. Yet by
Theorem~\ref{thm:infection}, $\#^\ell K$ is not $(n.5)$-solvable for
any $\ell\ge1$, since, by Lemma~\ref{signatureJ}, $\rho_\BZ (\#^N_{j=1}J)=N\rho_\BZ(J)$ is strictly more than $C_M$. Note that the Arf
invariant of $J$ is zero since our $4$-manifold $W$ is spin. Thus $\#^\ell K\notin \SF_{(n.5)}$ and it follows a fortiori that $\#^\ell K$ does not bound any
$(n+2.5)$-Grope in $D^4$.
\end{proof}

\section{Higher-Order Alexander Modules and Blanchfield Duality} \label{sec:noncommutative}

We have completed the proofs of our main theorems modulo the proof of Theorem~\ref{thm:injectivity}, which is a deep result concerning the nontriviality of the inclusion maps $\pi_1(M_R)\to \pi_1(W)$ where $W$ is an arbitrary ($n$)-solution for the knot $R$. To underscore the nature and difficulty of the task, note that it will be a direct consequence of Theorem~\ref{thm:injectivity} that, if $R$ is a genus two fibered ribbon knot and $W$ is the exterior in $B^4$ of \emph{any} slice disk for $R$, then $\pi_1(W)$ \emph{is not a solvable group}! For such a $W$ is an ($n$)-solution for every $n$ and if $\pi_1(W)^{(n)}=0$ then, for any $\eta_i\in \pi_1(M_R)^{(n)}$, $j_*(\eta_i)=0$. To prove Theorem~\ref{thm:injectivity} we will need the full power of the duality results of \cite[Section 2]{COT}, which we here review for the convenience of the reader.

In \cite{COT} it was observed that the classical Alexander module and the classical Blanchfield
linking form have higher-order generalizations that can
be defined using noncommutative algebra. These modules reflect the highly nonsolvable nature of the fundamental group of the knot exterior. Using these modules as measures of nonsolvability, it was shown that if a knot is a slice knot with slice disk $\Delta$, then $\pi_1(B^4-\Delta)$ is not arbitrary, but rather is constrained by $\pi_1(S^3-K)$. It is Poincar\'{e} Duality with twisted coefficients that provides connections between the higher-order modules of $B^4-\Delta$ and those of $S^3-K$. These constraints are expressed through higher-order Blanchfield duality. We also obtained similar results if $K$ is an ($n$)-solvable knot.

More specifically, recall that the classical Alexander module and Blanchfield form are associated to the infinite cyclic cover, $X_{\mathbb{Z}}$, of the exterior, $X$, of a knot $K$ in $S^3$, being merely $H_1(X_{\mathbb{Z}},\mathbb{Z})$ viewed as a module over $\mathbb{Z}[t^{\pm1}]$. More generally, consider any regular covering space, $X_{\G}$, that is obtained by taking \emph{iterated} covering spaces, each with torsion-free-abelian covering group. The group of covering translations, $\G$, of such a cover is \emph{poly-(torsion-free-abelian)} and hence solvable and torsion-free. We call $H_1(X_{\G},\mathbb{Z})$, viewed as a module over $\mathbb{Z}\G$, a \emph{higher-order Alexander module}. Even in the classical case, the ring $\mathbb{Z}[t^{\pm1}]$ can be simplified by localizing the coefficient group $\mathbb{Z}$ to get $\mathbb{Q}$ and considering the \emph{classical rational Alexander module}, $\SA_{\mathbb{Z}}$, a module over the PID $\mathbb{Q}[t^{\pm1}]$. Similarly, after appropriate localization of the ring $\mathbb{Z}\G$, the higher-order modules become finitely-generated torsion modules, $\SA_{\G}$ over a noncommutative PID $\mathbb{K}[t^{\pm1}]$ \cite[Proposition 2.11, Corollary 3.3]{COT}. Briefly, this is seen as follows. There is a split exact sequence 
\[
1\to\G '\to\G\to\mathbb{Z}\to 1
\]
 where $\G '$ is the commutator subgroup of $\G$. It follows that $\mathbb{Z}\G$ can be identified with a twisted Laurent polynomial ring $\mathbb{Z}\G'[t^{\pm1}]$. Because $\G$ is torsion-free and solvable, $\mathbb{Z}\G$ and $\mathbb{Z}\G'$ are right Ore domains and thus admit classical right quotient fields $\SK$ and $\mathbb{K}$. Let $\SR$ denote the ring obtained from $\mathbb{Z}\G$ by inverting the non-zero elements of $\mathbb{Z}\G'$, so that 
 \[
 \mathbb{Z}\G \subset \SR \subset \SK.
 \]
  The localized ring $\SR$ is then identified with the twisted Laurent polynomial ring $\mathbb{K}[t^{\pm1}]$ with coefficients in the skew field $\mathbb{K}$. Thus $\SR$ can seen to be a (noncommutative) PID by considering the obvious degree function on the polynomial ring. There is a classification theorem for modules over a noncommutative PID analogous to that for a commutative PID. Consequently the finitely generated modules 
  \[
  \SA_{\G}\cong \oplus _i (\SR/p_i(t)\SR)
  \]
then have a measurable ``size", namely their ranks over $\mathbb{K}$ (which are directly interpretable as the degrees of ``higher-order Alexander polynomials'' of the knot $K$!). These ranks are an important measure of nontriviality in the proofs of the next section. Henceforth we restrict our attention to these (\emph{localized}) higher-order Alexander modules.

Slightly more generally, if $X$ is any compact space and  $\phi:\pi_1(X)\to \G$ is any homomorphism where $\G$ is a poly-(torsion-free-abelian) group such that $\G/\G'\cong \mathbb{Z}$, then we may associate a (localized) higher-order Alexander module $\SA_{\G}(X)$ as above. Viewing $\phi$ as a system of twisted coefficients on $X$, note that $\SA_{\G}(X)$ is merely $H_1(X,\SR)\cong H_1(X,\mathbb{Z}\G)\otimes \SR$. We apply this to two spaces, first to $M$, the $3$-manifold obtained by $0$-framed surgery on $K$, and secondly to a $4$-manifold $W$ whose boundary is $M$. We see below that if $W$ is the complement of a slice disk for $K$ or more generally if $W$ is an $n$-solution for $M$ then the Alexander modules of $M$ and $W$ are strongly related.

Moreover we can define a higher-order linking form
$$ B\l_{\G}:\SA_{\G}(M):=H_1(M;\SR) \to  \Hom_{\SR}(\SA_{\G}(M),\SK/\SR) =: \SA_{\G}(M)^\# $$
as the composition of the Poincar\'e duality isomorphism to $H^2(M;\SR)$, the inverse of
the Bockstein to
$H^1(M;\SK/\SR)$, and the usual Kronecker evaluation map to $\SA_{\G}(M)^\#$. In the case that $\phi:\pi_1(M)\to \G$ is the abelianization map to $\mathbb{Z}$ we recover the classical $\SA_{\mathbb{Z}}$ and $B\l_{\mathbb{Z}}$.

\begin{thm}\label{higher Bforms}\cite[Theorem 2.13]{COT} There is a nonsingular hermitian linking form
$$ B\l_{\G}:H_1(M;\SR)\to H_1(M;\SR)^\#
$$ defined on the  higher-order Alexander module $\SA_{\G}(M)=H_1(M;\SR)$.
\end{thm}
\begin{thm}\label{self annihil} Suppose $M$ is 
$(n)$-solvable via $W$ and
$\phi:\pi_1(W)\to\G$ is a non-trivial coefficient system where
$\G^{(n)}=0$. Then the linking form
$B\l_{\G}(M)$ is hyperbolic, and in
fact the kernel $\SP$ of
\[
j_{\ast}:H_1(M;\SR)\to H_1(W;\SR)
\]
is self-annihilating ($\SP=\SP^{\bot}$ with respect to $B\l_{\G}$).
\end{thm}
It follows easily that the induced map $B\l_{\G}:\SP\to (\SA_{\G}(M)/\SP)^\#$ is an isomorphism \cite[Lemma 2.14]{COT}, and thus, by examining the exact sequence
$$
0\to\SP\to\SA_{\G}(M)\to \SA_{\G}(M)/\SP\to 0
$$
that $\rank_{\mathbb{K}}(\SP)=\rank_{\mathbb{K}}(\SA_{\G}(M)/\SP)=(1/2) \rank_{\mathbb{K}}(\SA_{\G}(M))$. Therefore we obtain the following crucial consequence of these theorems: the ``size'' of the \emph{image} of $j_{\ast}:H_1(M;\SR)\to H_1(W;\SR)$, as measured by the rank over $\mathbb{K}$, is \emph{precisely one half} of the ``size'' of $\SA_{\G}(M)$. In particular, as long as the latter is nontrivial, the image is nontrivial. But the ``size'' of the higher-order Alexander modules is easy to establish for fibered knots and turns out to be bounded below by $d-2$ where $d=\rank_{\mathbb{Q}}(\SA_{\mathbb{Z}}(M))$ is the degree of the classical Alexander polynomial of $K$ (whence the restriction in Theorem~\ref{thm:injectivity} to knots with $d=4$!).
In \cite[Corollary 4.8]{C}, this bound was indeed  established for {\em any} knot.

 The nontriviality that results from this sequence of ideas is the key element in the proof of our Theorem~\ref{thm:injectivity} in the next section, where more details will be given.

\section{Proof of Theorem~\ref{thm:injectivity}}\label{sec:proof}

Suppose $R$ is a genus 2 fibered knot in $S^3$ and $M$ is the result
of $0$-framed surgery on $R$. Note that $M$ fibers over $S^1$ with
fiber a closed genus 2 surface $\Sigma$. Let $S$ denote
$\1(\Sigma)$ which can be identified (under inclusion) with
$\1(M)\21$. Suppose also that $W$ is an
$(n)$-solution for $M$. The inclusion $j:M\lra W$ induces a map
$j:S\lra\1(W)\21$. Let $G=\1(W)\21=\1(W)\21_r$. From this
data, provided by the $(n)$-solution $W$, we shall abstract
certain algebraic properties and call this an {\bf
algebraic $(n)$-solution}. For the following definition we let
$F=F\<x_1,x_2,x_3,x_4\>$ be the free group and consider
$F\thra S$ in the standard way with $\{x_i\}$ corresponding to
a symplectic basis of $H_1(S)$ so that the kernel is normally
generated by $[x_1,x_2][x_3,x_4]$. We adopt the short-hand
$G_k=G/G\2k_r$.

\begin{defn}\label{algebraic} A homomorphism $r:S\to G$ is
called an {\bf algebraic $(n)$-solution} ($n\ge1$) if the
following hold:
\begin{enumerate}
\item $r_*:H_1(S;\BQ)\lra H_1(G;\mathbb{Q})$ has $2$-dimensional image
and after possibly re-ordering $\{x_1,x_2\}$ and $\{x_3,x_4\}$,
$r_*(x_1)$ and $r_*(x_3)$ are nontrivial.
 \item For each $0\le
k\le n-1$ the following composition is non-trivial even after
tensoring with the quotient field $\BK(G_k)$ of $\BZ G_k$:
$$
H_1(S;\BZ G_k)\overset{r_*}{\lra} H_1(G;\BZ G_k)\cong
G\2k_r/[G\2k_r,G\2k_r]\thra G\2k_r/G^{(k+1)}_r.
$$
\end{enumerate}
\end{defn}

We remark that if $r:S\lra G$ is an algebraic $(n)$-solution
then, for any $k<n$ it is an algebraic $k$-solution.

\begin{prop}\label{geomalg} The map $j:S\lra G$ (induced by
the $(n)$-solution $W$ above) is an algebraic $(n)$-solution.
\end{prop}

We postpone the proof of this proposition.

\begin{thm}\label{thm:paininbutt} For each $n\ge0$ there is a
finite set $\SP_n$ of pairs of elements of $F\2n$ with this
property: For any algebraic $(n)$-solution $r:S\lra G$ (no condition for $n=0$), at
least one such pair (which will be called a {\bf special
pair} for $r$) maps to a $\BZ G_n$-linearly independent set
under the composition:
$$
F\2n\lra S\2n/S\nn\cong H_1(S;\BZ S_n)\overset{r_*}{\lra}
H_1(S;\BZ G_n).
$$
\end{thm}

Assuming this theorem and the previous Proposition, we finish the
proof of Theorem~\ref{thm:injectivity}. Apply Theorem~\ref{thm:paininbutt}
to find a finite set $\SP_{n-1}$ of pairs of elements of
$F^{(n-1)}$. Since $S^{(n-1)}=\1(M)\2n$, the union of the elements
of $\SP_{n-1}$ is a finite set $\{\a_1,\dots,\a_m\}$ of elements
of $\1(M)\2n$ as required by Theorem~\ref{thm:injectivity}. Suppose
$W$ is an $(n)$-solution for $M$. Then by Proposition~\ref{geomalg}
the induced map $j:S\lra G$ is an algebraic $(n)$-solution. Suppose $n\geq 2$. Since
$j$ is also an algebraic $(n-1)$-solution, by
Theorem~\ref{thm:paininbutt}, at least one pair
$(y,z)\in\SP_{n-1}$ spans a $2$-dimensional subspace (over $\mathbb{K}(G_{n-1})$) of $H_1(S;\BK(G_{n-1}))$. But we claim that, if $n\geq 2$, $H_1(S;\BK(G_{n-1}))$ has rank $2$, so this subspace is all of $H_1(S;\BK(G_{n-1}))$. To establish this claim, observe that the Euler characteristic of $\Sigma$ can be computed using homology with coefficients in the skew field $\mathbb{K}(G_{n-1})$ (compare \cite[p.357]{C}). Thus $-2=\chi (\Sigma)=b_0-b_1+b_2$ where $b_i$ is the rank of $H_i(\Sigma;\BK(G_{n-1}))=H_i(S;\BK(G_{n-1}))$. The coefficient system $S\to G/G^{(1)}_r$ is nontrivial by property $(1)$ of the Definition~\ref{algebraic}. Hence $b_0=0$ (see \cite[Proposition 2.9]{COT}) and consequently $b_2$ is also zero by Poincar\'{e} Duality. Thus $b_1=2$ as claimed. Then, applying part 2 of the definition
of an algebraic $(n)$-solution (with $k=n-1$) we see that  at least
one of $\{y,z\}$ maps non-trivially under 
\[
F^{(n-1)}\lra S^{(n-1)}\overset{j}{\lra}G^{(n-1)}_r/G\2n_r.
\]
 If $n=1$ at least one of $\{x_1,x_2\}$ maps non-trivially by property $(1)$ of the definition of an algebraic $1$-solution. Hence, in any case, at
least one $\a_i$ has the property that $j_*(\a_i)\notin
G\2n_r=\1(W)^{\nn}_r$. Since each $\a_i$ actually comes from
$F^{(n-1)}$ and since $\1(M-\Sigma)\cong S$, we can represent the
$\a_i$ by simple closed curves in the complement of $\Sigma$, and
hence in the exterior of a Seifert surface (in $S^3$) for the knot
$R$. This is the collection
$\{\eta_i\}$ required by Theorem~\ref{thm:injectivity} whose proof is thus completed.

\begin{proof}[Proof of Proposition~\ref{geomalg}] Since $M$
is \ns\ via $W$, Theorem~\ref{thm:rank} below applies with $\G=\BZ$
and $n=1$ to show that $j_*:H_1(M;\BQ[t,t^{-1}])\lra
H_1(W;\BQ[t,t^{-1}])$ has rank $r/2$ over $\BQ$ where
$r=\rank_\BQ H_1(M_\infty;\BQ)$. But for a fibered knot,
$H_1(M;\BQ[t,t^{-1}])\cong H_1(M_\infty;\BQ)$ is equal to $H_1(S;\BQ)$ and, on the other hand,
$H_1(W;\BQ[t,t^{-1}])$ is given by 
\[
\1(W)\21/\1(W)\22\ox\BQ\cong
G/G\21\ox\BQ\cong H_1(G;\BQ).
\]
 For a genus $2$ fibered knot $r=4$
and so $j_*$ has $2$-dimensional image as required by condition 1.
Since the inclusion $M\overset{j}{\lra}W$ is an isomorphism on
$H_1$, there is a map $f:W\lra S^1$ such that $f^{-1}$ (regular
value) is an embedded $3$-manifold $Y$ whose boundary is $\Sigma$.
We have a factorization of $j_*$ as
$$
H_1(\Sigma;\BQ)\overset{i_*}{\lra}H_1(Y;\BQ)\overset{k_*}{\lra}
H_1(G;\BQ).
$$
 By the usual Poincar\'{e} Duality argument the kernel of $i_*$ is a Lagrangian of the intersection form on $H_1(\Sigma;\mathbb{Q)}$.
 In particular, in our setting, it has dimension $2$. This also implies that the image of $i_*$ has rank
2, so $\kernel j_*=\kernel i_*$. Suppose both $x_1$ and $x_2$ lay
in $\kernel j_*$ hence in $\kernel i_*$. This is a contradiction
since $x_1\cd x_2\neq0$ in $H_1(\Sigma;\BQ)$. For if $\<,\>$ is
the intersection form $H_2(Y,\p Y)\otimes H_1(Y)\lra\BQ$ and
$\p_*:H_2(Y,\p Y;\BQ)\lra H_1(\p Y;\BQ)$ then
$\<z,i_*x_2\>=\p_*z\cd x_2$ so if $x_1$ lay in $\kernel i_*$ then
$x_1$ would be of the form $\p_*z$ for some $z$, implying $x_1\cd
x_2=0$. Similarly at least one of $\{x_3,x_4\}$ has non-zero
image. This completes the verification of property 1.

Now suppose $k\le n-1$. Let $\G=\1(W)/\1(W)^{(k+1)}_r$ so $\G$
is PTFA and $\G\2n=\{e\}$. Then $\G ^{(1)}=G/G\2k_r\equiv G_k$.
By Theorem~\ref{thm:rank} below, the map $j_*:H_1(M;\BZ\G
^{(1)}[t^{\pm1}])\lra H_1(W;\BZ\G ^{(1)}[t^{\pm1}])$ has rank at
least 1 as a map of $\BZ\G ^{(1)}$-modules, i.e. $\BZ
G_k$-modules. But as $\BZ G_k$ modules this map is identical to
$j_*:H_1(M_\infty;\BZ G_k)\lra H_1(W_\infty;\BZ G_k)$ where
$W_\infty$ is the infinite cyclic cover of $W$. But
$H_1(M_\infty;\BZ G_k)=H_1(S;\BZ G_k)$ and $H_1(W_\infty;\BZ
G_k)=H_1(G;\BZ G_k)$ and we have shown that $j_*$ is non-trivial
even after tensoring with $\BK(G_k)$. Since the kernel of
$G\2k_r/[G\2k_r,G\2k_r]\lra G\2k_r/G^{(k+1)}_r$ is
$\BZ$-torsion, it is an isomorphism after tensoring with
$\BK(G_k)$ (which contains $\mathbb{Q})$. Thus $j$ is an algebraic
$(n)$-solution.
\end{proof}

\begin{thm}\label{thm:rank} Let $M$ be zero surgery on a knot.
Suppose $M$ is \ns\ via $W$ and $\psi:\1(W)\lra\G$ induces an
isomorphism upon abelianization where $\G$ is PTFA group and
$\G\2n=\{e\}$. Let $r=\rank_\BQ H_1(M_\infty;\BQ)$ where
$M_\infty$ is the infinite cyclic cover of $M$. Then
\[
j_*:H_1(M;\BK[t^{\pm1}])\lra H_1(W;\BK[t^{\pm1}])
\]
 has rank at least $(r-2)/2$ if $n>1$ and has rank $r/2$ if $n=1$, as a map of
$\BK$ vector spaces, where $\BK$ is the quotient field of $\BZ\G^{(1)}$.
\end{thm}

\begin{proof} Let $\SR=\BK[t^{\pm1}]$. By Theorem~\ref{higher Bforms} there exists a non-singular linking form
$\SB\ell:H_1(M;\SR)\lra H_1(M;\SR)^\#$. Let $\SA=H_1(M;\SR)$. By
Theorem~\ref{self annihil}, $P=\kernel(j_*)$ is an $\SR$-submodule
of $\SA$ which is self-annihilating with respect to $\SB\ell$.
It follows that the map $h:P\lra(\SA/P)^\#$ given by
$p\mapsto\SB\ell(p,\cd)$ is an isomorphism \cite[Lemma 2.14]{COT} . Note that $\SA/P$ is
isomorphic to the image of $j_*$. We claim that the rank over
$\BK$ of a finitely-generated right $\SR$-module $\SM$ is equal to
the $\BK$-rank of the right $\SR$-module
$\ov{\Hom_\SR(\SM,\SK/\SR)}\equiv\SM^\#$. Since $\SR$ is a
noncommutative PID \cite[Chapter 3]{J}\cite[Prop.4.5]{C}, any finitely-generated $\SR$-module is a
direct sum of cyclic modules \cite[Thm2.4 p.494]{Co}. Hence our claim can be seen by
examining the case of a cyclic module $\SM=\SR/p(t)\SR$ and
verifying that in this case $\SM^\#\cong(\ov{\SR/\SR
p(t)})\cong\SR/\bar p(t)\SR$ where $\bar p(t)$ is the result of
applying the involution from the group ring $\BZ\G$. One also
verifies that (just as in the commutative case) the $\BK$ rank of
such a cyclic module is the degree of the Laurent polynomial
$p(t)$. Since the degree of $\bar p(t)$ equals the degree of
$p(t)$, we are done. Hence $\rank_\BK(P)=\rank_\BK(\image j_*)$
and so this rank is at least $\rank_\BK(\SA)/2$. It remains only
to show that $\rank_\BK(\SA)$ is at least $r-2$ if $n>1$ and is
$r$ if $n=1$.

By hypothesis, $M$ is $0$-framed surgery on a knot $R$ in $S^3$.
Then $r$ is interpretable as the degree of the Alexander
polynomial of $R$. If $n=1$ then $\mathbb{K}=\mathbb{Q}$ and $\SA$
is the classical Alexander module, which is well known to have
$\BQ$-rank $r$. If $n>1$, by Corollary 4.8 of \cite{C}, $\rank_\BK
H_1(\sr;\BK[t^{\pm1}])\ge r-1$. Since $\SA$ depends only on
$\1(M)$ and the latter is obtained from $\1(\sr)$ by killing the
longitude, $\SA$ is obtained from $H_1(\sr;\BK[t^{\pm1}])$ by
killing the $\BK[t^{\pm1}]$-submodule generated by the longitude
$\ell$. If $n=1$ the longitude is trivial. Since $\ell$ commutes
with the meridian of $R$, 
\[
(t-1)_*[\ell]=0\in H_1(\sr;\BK[t^{\pm1}])
\]
 implying that this submodule is a
quotient of $\BK[t^{\pm1}]/(t-1)\BK[t^{\pm1}]\cong\BK$. Hence, if
$n>1$, $\rank_\BK\SA\ge r-2$.
\end{proof}

\noindent {\em Proof of Theorem~\ref{thm:paininbutt}}.  Set
$\SP_0=\{(x_1,x_2)\}$ and 
\[
\SP_1=\{([x_i,x_j],[x_i,x_k])| \
i,j,k \text{ distinct} \}. 
\]
Supposing $\SP_k$ has been defined,
define $\SP_{k+1}$ as follows. For each $(y,z)\in\SP_k$
include the following 12 pairs in $\SP_{k+1}$:
$([y,y^{x_i}],[z,z^{x_i}])$, $([y,z],[z,z^{x_i}])$,
$([y,y^{x_i}],[y,z])$ for $1\le i\le4$ and $y^x\equiv
x^{-1}yx$.

Now we fix $n$ and show $\SP_n$ satisfies the conditions of
the theorem. Fix an algebraic $(n)$-solution $r:S\lra G$. We
must show that there exists a special pair in $\SP_n$
corresponding to $r$. This is true for $n=0$ since
$H_1(S;\BZ G_0)\cong H_1(S;\BZ)$, so we assume $n\ge1$. Now
we need some preliminary definitions.

Let $F$ be the free group on $\{x_1,\dots,x_4\}$. Its
classifying space has a standard cell structure as a wedge of
4 circles $W$. Our convention is to consider its universal
cover $\wt W$ as a right $F$-space as follows. Choose a
preimage of the $0$-cell as basepoint denoted $*$. For each
element $w\in F\equiv\1(W)$, lift $w^{-1}$ to a path
$(\tl w^{-1})$ beginning at $*$. There is a unique deck
translation $\Phi(w)$ of $\wt W$ that sends $*$ to the
endpoint of this lift. Then $w$ acts on $\wt W$ by $\Phi(w)$.
This is the conjugate action of the usual left action as in
\cite{Ma}. Taking the induced cell structure on $\wt W$ and
tensoring with an arbitrary left $\BZ F$-module $A$ gives an
exact sequence
\begin{equation}\label{F}
0\lra H_1(F;A)\overset{d}{\lra} A^4\lra A\lra H_0(F;A)\lra0.
\end{equation}
Specifically consider $A=\BZ G$ where $\BZ F$ acts by left
multiplication via a homomorphism $\phi:F\lra G$. From the
interpretation of $H_1(F;\BZ G)$ as $H_1$ of a $G$-cover of $W$,
one sees that an element $g$ of $\ker(\phi)$ can be considered as
an element of $H_1(F;\BZ G)$. We claim that the composition
$\ker(\phi)\lra H_1(F;\BZ G)\overset{d}{\lra}(\BZ G)^4$ can be
calculated using the ``free differential calculus''
$\p=(\p_1,\dots,\p_4)$ where $\p_i:F\lra\BZ F$. Specifically we
assert that the diagram below commutes. Henceforth we abbreviate
maps of the form $(r,r,r,r):(\BZ F_n)^4\lra(\BZ G_n)^4$ by $r$.
\begin{equation}
\begin{diagram}\label{G}
\dgARROWLENGTH=1.2em
\node{F} \arrow[2]{e,t}{\partial} \node[2]{(\bbz F)^{4}}
\arrow{s,r}{\phi}\\
\node{ker(\phi)} \arrow{n} \arrow{e} \node{H_1(F;\bbz G)}
\arrow{e,t}{d} \node{(\bbz G)^{4}}
\end{diagram}
\end{equation}
where $\p_i(x_j)=\d_{ij}$, $\p_i(e)=0$ and
$\p_i(gh)=\p_ig+(\p_ih)g^{-1}$ for each $1\le i\le4$. This can be
seen as follows. Let $e_1$, $e_2$, $e_3$, $e_4$ denote lifts of
the 1-cells of $W$ to $\wt W$ that emanate from $*$ and are
oriented compatibly with $x_1,\dots,x_4$. For any word $g\in
F\equiv\1(W)$, $\p g$ is (by definition) the $\BZ F$-linear
combination of the $e_i$ that describes the 1-chain determined by
lifting a path representing $g$ to a path $\tl g$ in $\wt W$
starting at $*$. If $g\in\ker(\phi)$ then $g$ lifts to a loop in
the $\phi$-cover of $W$ (and $H_1(F;\BZ G)$ is the first homology
of this cover) and the 1-cycle it represents can be obtained by
pushing down $\tl g$ from the universal cover. In other words the
1-chain $\tl g$ in $C_1(\wt W)\ox_{\BZ F}\BZ G$ is obtained from
the 1-chain in $C_1(\wt W)$ by applying $\phi$ in each coordinate.
It only remains to justify the formula for $\p_i$. Note that the
usual formula for the standard left action is
\[
d_i(gh)=d_i(g)+gd(h).
\]
Our formula is obtained by setting
$\p_i=\bar d_i$ where $^-$ is the involution on the group ring.
Alternatively, this can be verified explicitly by induction on the
length of $h$. Suppose first that the length of $h$ is $1$. If
$h\neq x_i$ then the formula is clearly true, so consider
$\p_i(gx_i)$. The path $\tl g\tl x_i$, viewed as a 1-chain is
obviously equal to the 1-chain given by $\tl g$ (whose
$i^{\supth}$ coordinate is $\p_ig$) plus a certain translate of
$e_i$. Since the path from $*$ to the initial point of $\tl x_i$
is $\tl g$, this is the image of $e_i$ under the action of
$g^{-1}$ (in our convention). Hence $\p_i(gx_i)=\p_ig+g^{-1}$ as
claimed. Now it is a simple matter to verify the inductive step by
expressing $h=h'\cd h''$ where $h'$ and $h''$ are of lesser
length. This is left to the reader.

Note that $r\pi_k:F\to F_k\to G_k$
is the same as $\pi_k r:F\to S\to G\to G_k$.
\begin{lem}\label{finally} Given an algebraic $(n)$-solution ($n\geq 1$)
$r:S\lra G$, after reordering $\{x_1,x_2\}$ and $\{x_3,x_4\}$ so that $r_*(x_1)$ and $r_*(x_3)$ are non-trivial as per part $(1)$ of Definition~\ref{algebraic}, for each $k$, $1\le k\le n$, there is at least one pair
$(y,z)\in\SP_k$ with the following {\bf good} properties:
\begin{enumerate}
\item $\p_4y=\p_4z=0$ \item The vectors
$(r\pi_k\p_2y,r\pi_k\p_3y)$ and $(r\pi_k\p_2z,r\pi_k\p_3z)$ (i.e.
the vectors consisting of the second and third
coordinates of the images of $y$ and $z$ under the composition
$F\2k\overset{\pi_k\p}{\lra}(\BZ F_k)^4\overset{r}{\lra}(\BZ
G_k)^4)$ are right linearly independent over $\BZ G_k$. (Note that
property $(1)$ ensures that the fourth coordinates are
zero).
\end{enumerate}
\end{lem}
\begin{proof}[Proof that Lemma~\ref{finally} $\Longrightarrow$
Theorem~\ref{thm:paininbutt}] The set $\SP_n$ was defined above.
Given an algebraic $(n)$-solution $r:S\lra G$, re-order $\{x_1,x_2\}$ and $\{x_3,x_4\}$ so that $r_*(x_1)$ and $r_*(x_3)$ are non-trivial as is possible by part $(1)$ of Definition~\ref{algebraic}. Then the Lemma provides a pair
$(y,z)\in\SP_n$ that has the listed {\bf good} properties with respect to $r$. We
verify that $(y,z)$ is a {\bf special pair} with respect to $r$. Consider the
diagram below. Recall $G_n=G/G\2n_r$.

\begin{equation} \begin{diagram}
\node{F\2n} \arrow{e}\arrow[2]{s,r}{r\2n} \node{H_1(F;\BZ F_n)} \arrow{e,t}{d} \arrow{s,r}{r_*} \node{(\BZ F_n)^4} \arrow{s,r}{(r_n)^4}\\
\node[2]{H_1(F;\BZ G_n)} \arrow{e,t}{d'} \arrow{s,r}{i_*} \node{(\BZ G_n)^4} \arrow{s}\\
\node{S\2n}\arrow{e} \node{H_1(S;\BZ G_n)} \arrow{e,t}{d''} \node{ (\BZ
G_n)^4/\<d'(T)\>} 
\end{diagram} \end{equation}

The horizontal composition on top is $\pi_n\p$. The right-top
square commutes by naturality of the sequence (\ref{F}) above. Let $T\in
F$ denote the single relation such that $S=F/\<T\>$. Since $T$ is
in the kernel of $F\lra S\lra G_n$, it represents an element of
$H_1(F;\BZ G_n)$ and it generates the kernel of the epimorphism
$i_*$. Since $d$ and $d'$ are monomorphisms, $d''$ is a
monomorphism. Hence to show $(y,z)$ is special, it suffices to
show that the set of three $4$-tuples $\{r\pi_n\p(y),r\pi_n\p(z),d'T\}$ is $\BZ G_n$-linearly
independent in $(\BZ G_n)^4$. This will follow immediately from the {\bf good}
properties of $y$ and $z$ once we verify that the $4^{\supth}$
coordinate of $d'T$ is non-zero. The standard relation $T$ after our possible re-ordering, is either $g[x_3,x_4]$ or $g[x_4,x_3]$ where $g$ is either $[x_1,x_2]$ or $[x_2,x_1]$. If $d$ stands for any one of the $\p_i$ then one calculates that $d(g^{-1})=-(dg)g$ and 
\[
d([g,h])=dg+(dh)g^{-1}-(dg)gh^{-1}g^{-1}-(dh)hgh^{-1}g^{-1}.
\]
 Using these one calculates that $\p_4(g[x_3,x_4])=(x_3^{-1}-[x_4,x_3])g^{-1}$ and 
$\p_4(g[x_4,x_3])=(1-x_4x^{-1}_3x^{-1}_4)g^{-1}$. The
$4^{\text{th}}$ coordinate of $d'T$ is $r\pi_n\p_4(T)$ by
diagram (\ref{G}). If this vanished in $\BZ G_n$ then its image $r\pi_1\p_4(T)$ would certainly vanish in
$\BZ G_1=\BZ[G/G\21_r]$ (recall that $n\geq 1$). But $r\pi_1\p_4(T)$ is either $r(x_3^{-1})-1$ or $1-r(x_3^{-1})$, which can only vanish if $r(x_3)$ is
trivial in $G/G\21_r=H_1(G;\BZ)/$torsion, i.e. $r_*(x_3)=0$ in
$H_1(G;\BQ)$, contradicting our choice of $x_3$. Thus the $4^{\supth}$ coordinate of $d'T$ is non-trivial. This
completes the verification that the Lemma implies the Theorem.
\end{proof}

\begin{proof}[Proof of Lemma~\ref{finally}] The integer $n\geq 1$
is fixed throughout. By definition of an algebraic $(n)$-solution we may re-order so that
$r_*(x_1)$ and $r_*(x_3)$ are non-trivial in
$H_1(G;\BQ)=G_1\ox\BQ$. This implies that $r\1(x_1)$ and
$r\1(x_3)$ are non-trivial in $G_1$. We will prove the Lemma by
induction on $k$. We begin with $k=1$. Consider the pair
\[
(y,z)=([x_1,x_2],[x_1,x_3])\in\SP_1.
\]
 We claim that $(y,z)$ has
the {\bf good} properties. Certainly $\p_4y=\p_4z=0$ since $x_4$
does not appear in the words $y$ and $z$. Similarly the
third coordinate of the image of $y$ and the second
coordinate of the image of $z$ are zero. Hence to establish the
second good property, it suffices to show that $r\1\p_2y\neq0\neq
r\1\p_3z$. Since $\p_2y=x^{-1}_1-[x_2,x_1]$ and
$\p_3z=x^{-1}_1-[x_3,x_1]$ this follows since $r\1(x^{-1}_1)$ and
$r\1(x^{-1}_3)$ are non-trivial in $G_1$. Therefore the base of the
induction $(k=1)$ is established.

Now suppose the conclusions of the Lemma have been established for
$1,\dots,k$ where $k<n$. We establish them for $k+1$. Let
$(y,z)\in\SP_k$ be a pair that has the good properties. This means
that $\p_4y=\p_4z=0$ and that the vectors
$(r\pi_k\p_2y,r\pi_k\p_3y)$ and $(r\pi_k\p_2z,r\pi_k\p_3z)$ are
(right) linearly independent over $\BZ G_k$. Consider the
following 3 elements of $\SP_{k+1}$: 
\[
([y,y^{x_1}],[z,z^{x_1}]),\quad ([y,z],[z,z^{x_1}]), \quad ([y,y^{x_1}],[y,z])
\]
where $y^x=x^{-1}yx$.
We will show that at least one of these pairs $(y_{k+1},z_{k+1})$
has the good properties, finishing the inductive proof of
Lemma~\ref{finally}.

First note that in all cases $\p_4y_{k+1}=\p_4z_{k+1}=0$ since
$\p_4y=\p_4z=0$. For the remainder of this proof we
write $x$ for $x_1$, suppressing the subscript. We need to show
that there is at least one of the 3 pairs $(y_{k+1},z_{k+1})$ such
that the vectors $(r\pi_{k+1}\p_2y_{k+1},r\pi_{k+1}\p_3y_{k+1})$
and $(r\pi_{k+1}\p_2z_{k+1},r\pi_{k+1}\p_3z_{k+1})$ are $\BZ
G_{k+1}$-linearly independent.

\noindent{\bf Case 1}: {\sl Both $r\pi_{k+1}(y)$ and
$r\pi_{k+1}(z)$ are non-trivial in $G_{k+1}$}.

In this case we will show that the pair
$(y_{k+1},z_{k+1})=([y,y^x],[z,z^x])$ satisfies property $(2)$. Let
$d$ be either $\p_2$ or $\p_3$. One has $dx=0$ and $dy^x=(dy)x$. Using this and our previous computations of $d([g,h])$, one computes that $d([y,y^x])=(dy)p$
where 
\[
p=1+xy^{-1}-(y^x)^{-1}[y^x,y]-x[y^x,y].
\]
 Similarly $d([z,z^x])=(dz)q$ where $q=1+xz^{-1}-(z^x)^{-1}[z^x,z]-x[z^x,z]$.
We must show the vectors $(r\pi_{k+1}((\p_2y)p)$,
$r\pi_{k+1}((\p_3y)p))$ and $(r\pi_{k+1}((\p_2z)q)$,
$r\pi_{k+1}((\p_3z)q))$ are $\BZ G_{k+1}$-linearly independent.
Note that the first vector is a right multiple of
$v_{k+1}=(r\pi_{k+1}\p_2y$, $r\pi_{k+1}\p_3y)$ by $r\pi_{k+1}p$
and the second is a right multiple of
$w_{k+1}=(r\pi_{k+1}\p_2z$, $r\pi_{k+1}\p_3z)$ by $r\pi_{k+1}q$. The
right factor $r\pi_{k+1}p$ is seen to be non-trivial in $\BZ
G_{k+1}$ as follows. First observe that $[y,y^x]\in F^{(k+1)}$ so $r\pi_{k+1}([y,y^x])=1$. Then note
\[
r\pi_{k+1}p=r\pi_{k+1}(1+xy^{-1}-(y^x)^{-1}-x)
\]
is a linear
combination of 4 group elements $e$, $r\pi_{k+1}(xy^{-1})$,
$r\pi_{k+1}((y^x)^{-1})$, and $r\pi_{k+1}(x)$ in $G_{k+1}$. For
$r\pi_{k+1}p$ to vanish in $\BZ G_{k+1}$ the elements must pair up
in a precise way and in particular such that
$r\pi_{k+1}(x)=r\pi_{k+1}(xy^{-1})$ in $G_{k+1}$. This is a
contradiction since $r\pi_{k+1}(y)\neq e$ by hypothesis. No other
pairing is possible because the projections of the 4 elements to
$G_1$ are $e$, $r\1(x)$, $e$ and $r\1(x)$ and we have already
noted that $r\1(x)$ is non-trivial in $G_1$. An entirely
similar argument shows that the right factor $r\pi_{k+1}q$ is
non-trivial, using the non-triviality of $r\pi_{k+1}(z)$. Since
these right factors are non-trivial and $\BZ G_{k+1}$ has no zero
divisors, the linear independence of $\{v_{k+1},w_{k+1}\}$ would
be sufficient to imply the linear independence of the original set
of vectors. Recall that our hypothesis on $(y,z)$ ensures that the
set 
\[
\{v_k,w_k\}=\{(r\pi_k\p_2y,r\pi_k\p_3y), \ (r\pi_k\p_2z,r\pi_k\p_3z)\}
\]
 is linearly independent in $\BZ
G_k$. Note that $v_k$ and $w_k$ are the images of $v_{k+1}$ and
$w_{k+1}$ under the canonical projection $(\BZ G_{k+1})^2\lra(\BZ
G_k)^2$. We assert that the linear independence of $\{v_k,w_k\}$
implies the linear independence of $\{v_{k+1},w_{k+1}\}$ since the
kernel of $G_{k+1}\lra G_k$ is a torsion free abelian group.
Details follow. This will complete the verification, under the
assumptions of Case 1, that at least one pair $(y_{k+1},z_{k+1})$
has {\bf good} properties.

To establish our assertion above, consider more generally an
arbitrary endomorphism $f$ of free right $\BZ G_{k+1}$ modules
$f:(\BZ G_{k+1})^2\lra(\BZ G_{k+1})^2$ given by 
$\begin{pmatrix} a
&b\\ c &d\end{pmatrix} 
\begin{pmatrix} x\\ y\end{pmatrix}=(ax+by,cx+dy)$ for
$a$, $b$, $c$, $d$, $x$, $y\in\BZ G_{k+1}$. Then $f$ induces $\bar f:(\BZ G_k)^2\lra(\BZ G_k)^2$ given by the matrix
$\begin{pmatrix} \bar a  &\bar b\\ \bar c  &\bar d\end{pmatrix}$
where $\bar a$ is the projection of $a$. Suppose that $\bar f$ is
injective. Let $H=G\2k_r/G^{(k+1)}_r$ and note that $\BZ
G_{k+1}\ox_{\BZ[H]}\BZ\cong\BZ G_k$ where $a\ox1\mapsto\bar a$, as
$\BZ G_{k+1}-\BZ$ bimodules. Moreover (under this identification)
$f$ descends to 
\[
f\ox\id:(\BZ G_k)^2\lra(\BZ G_k)^2
\]
 sending $(\bar z,\bar w)$ (for $z$, $w\in\BZ G_{k+1}$) to ($\bar a\bar
z+\bar c\bar w$, $\bar b\bar z+\bar d\bar w$), thus agreeing with
$\bar f$ above. Since $H$ is torsion-free-abelian, a theorem of
Strebel \cite[Section 1]{Str} ensures that the injectivity of
$f\ox\id=\bar f$ implies the injectivity of $f$. An application of
this general fact with $\{(a,c),(b,d)\}=\{v_{k+1},w_{k+1}\}$ shows
the latter set is right linearly independent.

\noindent{\bf Case 2}: $r\pi_{k+1}(y)=e$ {\sl and
$r\pi_{k+1}(z)\neq e$ in $G_{k+1}$}.

In this case we claim that the pair
$(y_{k+1},z_{k+1})=([y,z],[z,z^x])$ satisfies the {\bf good}
property $(2)$. We see that
$dy_{k+1}=dy(1-z^{-1}[z,y])+dz(y^{-1}-[z,y])$ where $d=\p_2$ or
$d=\p_3$. Thus
\[
r\pi_{k+1}dy_{k+1}=(r\pi_{k+1}dy)(1-r\pi_{k+1}(z^{-1})).
\]
Therefore, one of the vectors,
$(r\pi_{k+1}\p_2y_{k+1},r\pi_{k+1}\p_3y_{k+1})$ is a right
multiple of the vector $v_{k+1}=(r\pi_{k+1}\p_2y,r\pi_{k+1}\p_3y)$
by a non-zero divisor $1-r\pi_{k+1}(z^{-1})$. Hence, just as in
Case 1, we can abandon the former and consider $v_{k+1}$. Recall
also that $dz_{k+1}=(dz)q$ as in Case 1. One checks that
$r\pi_{k+1}q\neq0$ in $\BZ G_{k+1}$ just as in Case 1 using the
fact that $r\pi_{k+1}(z)\neq e$. Thus we can reduce to considering
the vector $w_{k+1}=(r\pi_{k+1}\p_2z,r\pi_{k+1}\p_3z)$ as above.
We finish the proof of Case 2 just as in Case 1, using the fact
that our hypothesis guarantees that the vectors called $w_k$ and
$v_k$ are linearly independent in $\BZ G_k$.

\noindent{\bf Case 3}: $r\pi_{k+1}(y)\neq e$ {\sl and
$r\pi_{k+1}(z)=e$ in $G_{k+1}$}.

In this case we claim that the pair
$(y_{k+1},z_{k+1})=([y,y^x],[y,z])$ satisfies the {\bf good}
property $(2)$. Following Case 2,
$r\pi_{k+1}dz_{k+1}$ is equal to $(r\pi_{k+1}dz)(r\pi_{k+1}(y^{-1})-1)$ where
$d=\p_2$ or $\p_3$. Moreover $dy_{k+1}=(dy)p$. One finishes just
as in Case 2.

\noindent{\bf Case 4}: $r\pi_{k+1}(y)=r\pi_{k+1}(z)=e$ {\sl
in $G_{k+1}$}.

We claim that this case cannot occur. For recall that, by the inductive hypothesis, the pair
$(y,z)\in P_k$ has the {\bf good} properties for the algebraic $(n)$-solution $r$ where $1\le k\le
n-1$. But $r$ is also an algebraic $k$-solution since $k<n$ and by the proof of Lemma~\ref{finally}
$\Longrightarrow$ Theorem~\ref{thm:paininbutt}, the pair $(y,z)$ is a \textbf{special pair} for $r$. Thus under the composition
$$
F\2k\lra S\2k\lra H_1(S;\BZ S_k)\overset{r_*}{\lra}H_1(S;\BZ
G_k)
$$
the set $\{y,z\}$ maps to a linearly independent set. Since
$H_1(S;\BZ G_k)$ has rank 2 (as we showed in the paragraph following the statement of Theorem~\ref{thm:paininbutt}), this set is a basis of
$H_1(S;\BK(G_k))$. Since $r$ is also an algebraic
$(n)$-solution and $k\leq n-1$, by property $(2)$ of Definition~\ref{algebraic} the composition of the above
with the map
$$
H_1(S;\BZ G_k)\overset{r_*}{\lra}H_1(G;\BZ G_k)\cong
G\2k_r/G^{(k+1)}_r
$$
is non-trivial when restricted to $\{y,z\}$. On the other
hand this combined map $F\2k\lra G\2k_r/G^{(k+1)}_r$ is
clearly given by 
\[
y\mapsto r\pi_{k+1}(y) \text{ and } z\mapsto r\pi_{k+1}(y)
\]
 so it is not possible that both $r\pi_{k+1}(y)$
and $r\pi_{k+1}(z)$ lie in $G^{(k+1)}_r$. Therefore Case 4
is not possible.

This completes the proof that one of the 3 new pairs
$(y_{k+1},z_{k+1})$ satisfies the {\bf good} properties and
hence concludes the inductive step of our proof of Lemma~\ref{finally}, as well as Theorem~\ref{thm:paininbutt}.
\end{proof}


\end{document}